\documentclass[aop,MSNbibl,dvips]{arximspdf}
\usepackage{mathbh}
\usepackage{graphicx}

\doi{10.1214/12-AOP799} 
\volume{42}
\issue{2}
\pubyear{2014}
\firstpage{725}
\lastpage{759}

\makeatletter

\newcommand{\rrvert}{\vert}
\newcommand{\llvert}{\vert}
\def\textnormal{\mathrm}

\newcommand{\eqref}[1]{(\ref{#1})}
\newtheorem{thmm}{Theorem}[section]
\newtheorem{lem}[thmm]{Lemma}
\newtheorem{prop}[thmm]{Proposition}

\newproclaim{defn}[thmm]{Definition}
\newproclaim{ex}[thmm]{Example}
\newproclaim{conj}[thmm]{Conjecture}
\newproclaim{rem}[thmm]{Remark}
\newproclaim{rems}[thmm]{Remarks}
\newproclaim{ques}[thmm]{Question}

\def\L{\mathbb{L}}

\def\P{\mathbb{P}}
\def\Pmu{\mathbb{P}_\mu}

\def\N{\mathbb{N}}
\def\E{\mathbb{E}}
\def\D{\mathbb{D}}
\def\Db{\overline{\mathbb{D}}}
\def\C{\mathcal{C}}
\def\R{\mathbb{R}}
\def\Q{\mathbb{Q}}
\def\D{\mathbb{D}}
\def\T{\mathbb{T}}
\def\F{\mathcal{F}}

\def\Ln{\Lambda_{\t_n}}
\def \t {\mathfrak{t}}
\def\card{\operatorname{Card}}

\def\exc{\mathrm{exc}}

\def\cl{\operatorname{cl}}
\def\a{\alpha}
\def\b{\beta}
\def\me{\mathbh{e}}

\def\z{\zeta}

\makeatother

\begin{document}
\begin{frontmatter}

\title{Random stable laminations of the disk}
\runtitle{Random stable laminations}

\begin{aug}
\author{\fnms{Igor} \snm{Kortchemski}\corref{}\ead[label=e1]{igor.kortchemski@normalesup.org}}
\runauthor{I. Kortchemski}
\affiliation{Universit\'{e} Paris-Sud}
\address{Laboratoire de math\'{e}matiques\\
UMR 8628 CNRS\\
Universit\'{e} Paris-Sud\\
91405 ORSAY Cedex\\
France\\
\printead{e1}}
\end{aug}

\received{\smonth{2} \syear{2012}}
\revised{\smonth{8} \syear{2012}}

%
\begin{abstract}We study large random dissections of polygons. We
consider random
dissections of a regular polygon with $n$ sides, which are chosen
according to Boltzmann weights in the domain of attraction of a
stable law of index $\theta\in(1,2]$. As $n$ goes to infinity, we
prove that these random dissections converge in distribution toward
a random compact set, called the random stable lamination. If
$\theta=2$, we recover Aldous' Brownian triangulation. However, if
$\theta\in(1,2)$, large faces remain in the limit and a different
random compact set appears. We show that the random stable lamination can
be coded by the continuous-time height function associated to the
normalized excursion of a strictly stable spectrally positive L\'{e}vy
process of index $\theta$. Using this coding, we establish that the
Hausdorff dimension of the stable random lamination is almost surely
$2-1/\theta$.
\end{abstract}

%


\end{frontmatter}

\section*{Introduction}

In this article we study large random dissections of polygons. A~\emph{dissection} of a polygon is the union of the sides of the
polygon and of a collection of diagonals that may intersect only at
their endpoints. The faces are the connected components of the
complement of the dissection in the polygon. The particular case of
triangulations (when all
faces are triangles) has been extensively studied in the literature.
For every integer $n \geq3$, let $P_n$ be the regular polygon with~$n$
sides whose vertices are the $n$th roots of unity. It is well
known that the number of triangulations of $P_n$ is the Catalan
number of order $n-2$. In the general case, where faces of degree
greater than three are allowed, there is no known explicit formula
for the number of dissections of $P_n$, although an asymptotic
estimate is known (see \cite{FlajoletNoy,CK}). Probabilistic aspects
of uniformly distributed random triangulations have been
investigated; see, for example, the articles \cite{GaoWormald1,GaoWormald2}
which study graph-theoretical properties of uniform triangulations
(such as the maximal vertex degree or the number of vertices of
degree $k$). Graph-theoretical properties of uniform dissections of
$P_n$ have also been studied, extending the previously mentioned
results for triangulations (see \cite{Berna,CK}).

From a more geometrical point of view, Aldous
\cite{Aldous,Aldous2} studied the shape of a large uniform triangulation
viewed as a random compact subset of the closed unit disk. See also the
work of
Curien and Le Gall \cite{CLG}, who discuss a random continuous
triangulation (different from Aldous' one) obtained as a limit of
random dissections constructed recursively. Our goal is to
generalize Aldous' result by studying the shape of large random
dissections of $P_n$, viewed as random variables with values in the
space of all compact subsets of the disk, which is equipped with the
usual Hausdorff metric.


Let us state more precisely Aldous' results. Denote by $\mathfrak
{t}_n$ a
uniformly distributed random triangulation of $P_n$.
There exists a random compact subset $\mathfrak{t}$ of the closed unit disk
$\Db$ such that the sequence $(\mathfrak{t}_n)$ converges in distribution
toward $\mathfrak{t}$. The random compact set $\mathfrak{t}$ is a continuous
triangulation, in the sense that $\Db\setminus\mathfrak{t}$ is a disjoint
union of open triangles whose vertices belong to the unit circle.
Aldous also explains how $\mathfrak{t}$ can be explicitly constructed using
the Brownian excursion and computes the Hausdorff dimension of
$\mathfrak{t}$,
which is equal almost surely to $3/2$ (see also \cite{LGP}).

In this work, we propose to study the following generalization of
this model. Consider a probability distribution $(\mu_j)_{j \geq0}$
on the nonnegative integers such that $\mu_1=0$ and the mean of $ \mu$
is equal to $1$. We suppose that $\mu$ is in
the domain of attraction of a stable law of index $\theta\in
(1,2]$. For every integer $n\geq2$, let $\L_n$ be the set of all
dissections of $P_{n+1}$, and consider the following Boltzmann
probability measure on $\L_n$ associated to the weights~$(\mu_j)$:
\[
\mathbb{P}^ {\mu}_ {n}(\omega)=\frac{1}{Z_n} \prod
_{f\ \mathrm{face\ of
}\ \omega} \mu_{\deg(f)-1},\qquad \omega\in
\L_n,
\]
where $\deg(f)$ is the
degree of the face $f$, that is, the number of edges in the boundary of
$f$, and $Z_n$ is a normalizing constant. Note that the definition of
$\mathbb{P}^ {\mu}_ {n}$ involves only $\mu_2, \mu_3, \ldots,$
and $
\mu_0$ is the missing constant to obtain a probability measure. Under
appropriate conditions on $ \mu$, this definition makes sense for all
sufficiently large integers $n$. Let us mention two important special
cases. If $ \mu_0 = 2 - \sqrt{2}$ and $ \mu_i= ((2 - \sqrt{2})/2)^
{i-1}$ for every $i \geq2$, one easily checks that $ \mathbb{P}^
{\mu
}_ {n}$ is uniform over $ \L_n$. If $ p \geq3$ is an integer and if $
\mu_0 = 1-1/(p-1)$, $ \mu_ {p-1}=1/(p-1)$ and $ \mu_i=0$ otherwise,
$\mathbb{P}^ {\mu}_ {n}$ is uniform over dissections of $ \L_n$ with
all faces of degree $p$ (in that case, we must restrict our attention
to values of $n$ such that $n-1$ is a multiple of $p-2$, but our
results carry over to this setting).

We are
interested in the following problem. Let $\mathfrak l_n$ be a random
dissection distributed according to $\mathbb{P}_n$. Does the
sequence $(\mathfrak l_n)$ converge in distribution to a random compact
subset of $\Db$? Let us mention that this setting is inspired by
\cite{LGM}, where Le Gall and Miermont consider random planar maps
chosen according to a Boltzmann probability measure, and show that if the
Boltzmann weights do not decrease sufficiently fast, large faces
remain in the scaling limit. We will see that this phenomenon occurs
in our case as well.

In our main result Theorem \ref{thmcv1}, we first consider the
case where the variance
of $\mu$ is finite and then show that $\mathfrak l_n$ converges
in distribution to Aldous' Brownian
triangulation as $n \rightarrow\infty$. This extends Aldous'
theorem to random dissections which are not necessarily
triangulations. For instance, we may let $\mathfrak l_n$ be uniformly
distributed over the set of all dissections whose faces are all quadrangles
(or pentagons, or hexagons, etc.). As noted above, this requires that
we restrict our attention
to a subset of values of $n$, but the convergence of $\mathfrak l_n$ toward
the Brownian triangulation still holds. This maybe surprising result
comes from the fact that certain sides of the squares (or of the
pentagons, or
of the hexagons, etc.) degenerate in the limit. See also the recent
paper \cite{CK} for other classes of noncrossing configurations of the
polygon that converge to the Brownian triangulation.

On the other hand, if $\mu$ is in the domain of attraction of
a stable law of index $\theta\in(1,2)$, Theorem \ref{thmcv1} shows
that $(\mathfrak l_n)$
converges in distribution to another random compact subset $\mathfrak l$
of $\Db$, which we call the $\theta$-stable random lamination of the
disk. The random compact subset $\mathfrak l$ is the union of the unit circle
and of infinitely many
noncrossing chords, which can be constructed as follows. Let
$X^{\mathrm{exc}}=(X^{\mathrm{exc}}_t)_{0 \leq t \leq1}$ be the
normalized excursion of the
strictly stable spectrally positive L\'{e}vy process of index $\theta$
(see Section \ref{sec2.1} for a precise definition). For $0 \leq s < t \leq
1$, we set $ s \mathop{\simeq}^{X^{\mathrm{exc}}} t$ if $t= \inf\{ u
> s; X^{\mathrm{exc}}_u \leq X^{\mathrm{exc}}_{s-}\}$, and $s
\mathop{\simeq}^{X^{\mathrm{exc}}} s$ by convention. Then
%
\begin{equation}
\label{eqreprX}\mathfrak l= \bigcup_{s \mathop{\simeq}^{X^{\mathrm{exc}}} t}
\bigl[e^{-2 \textnormal
{i} \pi s},e^{-2
\textnormal{i} \pi
t}\bigr],
\end{equation}
where $[u,v]$ stands for the line segment between
the two complex numbers $u$ and~$v$. In particular, the latter set
is compact, which is not obvious {a priori}.

In order to study fine properties of the set $\mathfrak l$, we derive an
alternative representation in terms of the so-called height process
$H^{\mathrm{exc}}=(H^{\mathrm{exc}}_t)_{0 \leq t \leq1}$ associated
with $X^{\mathrm{exc}}$ (see
\cite{Duquesne,DuquesneLG} for the definition and properties of
$H^{\mathrm{exc}}$). Note that $H^{\mathrm{exc}}$ is a random
continuous function on $[0,1]$
that vanishes at $0$ and at $1$ and takes positive values on
$(0,1)$. Then Theorem \ref{thmXH} states that
%
\begin{equation}
\label{eqreprH}\mathfrak l= \bigcup_{s
\mathop{\thickapprox}^{H^{\mathrm{exc}}} t}
\bigl[e^{-2 \textnormal{i}
\pi s},e^{-2
\textnormal{i} \pi
t}\bigr],
\end{equation}
where, for $s,t \in[0,1]$, $s
\mathop{\thickapprox}^{H^{\mathrm{exc}}} t$ if $H^{\mathrm
{exc}}_s=H^{\mathrm{exc}}_t$ and $H^{\mathrm{exc}}_r > H^{\mathrm
{exc}}_s$ for
every $r \in( s \wedge t, s \vee t)$, or if $(s,t)$ is a limit of
pairs satisfying these properties. This is very closely related to
the equivalence relation used to define the so-called stable tree,
which is coded by $H^{\mathrm{exc}}$ (see
\cite{Duquesne}). The representation \eqref{eqreprH} thus shows that
the $\theta$-stable random lamination is connected to the $\theta$-stable
tree in the same way as the Brownian triangulation is connected
to the Brownian CRT (see \cite{Aldous2} for applications of the latter
connection).
The representation \eqref{eqreprH} also allows us to establish that
the Hausdorff dimension of $\mathfrak l$ is almost surely equal to
$2-1/\theta$. Note that for $\theta=2$, we obtain a Hausdorff
dimension equal to $3/2$, which is consistent with Aldous' result.
Additionally, we verify that the Hausdorff dimension of the set of endpoints
of all chords in $\mathfrak l$ is equal to $1-1/\theta$.

Finally, we derive precise information about the faces of
$\mathfrak l$, which are the connected components of the complement
of $\mathfrak l$ in the closed unit disk. When $\theta=2$, we already
noted that
all faces
are triangles. On the other hand, when $\theta\in(1,2)$, each face is
bounded by infinitely many chords. We prove more precisely that
the set of extreme points of the closure of a face (or, equivalently, the
set of points of the closure that lie on the circle) has Hausdorff
dimension $1/\theta$.

Let us now sketch the main techniques and arguments used to
establish the previous assertions. A key ingredient is the fact that
the dual graph of $\mathfrak l_n$ is a Galton--Watson tree conditioned on
having $n$ leaves. In our previous work \cite{K}, we establish limit
theorems for Galton--Watson trees conditioned on their number of
leaves and, in particular, we prove an invariance principle
stating that the rescaled Lukasiewicz path of a Galton--Watson tree
conditioned on having $n$ leaves converges in distribution to
$X^{\mathrm{exc}}$
(see Theorem \ref{thmleaves} below). Using this result, we are able
to show that $\mathfrak l_n$ converges toward the random compact set
$\mathfrak l$
described by~(\ref{eqreprX}). The representation \eqref{eqreprH}
then follows from relations between $X^{\mathrm{exc}}$ and $H^{\mathrm
{exc}}$. Finally, we use
\eqref{eqreprH} to verify that the Hausdorff dimension of $\mathfrak
l$ is
almost surely equal to $2-1/\theta$. This calculation relies in part
on the time-reversibility of the process $H^{\mathrm{exc}}$. It seems more
difficult to derive the Hausdorff dimension of $\mathfrak l$ from the
representation \eqref{eqreprX}.

The paper \cite{CK} develops a number of applications of the present
work to enumeration problems and asymptotic properties of uniformly
distributed random dissections.

The paper is organized as follows. In Section \ref{sec1} we present the
discrete framework. In particular, we introduce Galton--Watson trees
and their coding functions. In Section \ref{sec2} we discuss the normalized
excursion of the strictly stable spectrally positive L\'{e}vy process of
index $\theta$ and its associated lamination $L(X^{\mathrm{exc}})$.
In Section~\ref{sec3}
we prove that $(\mathfrak l_n)$ converges in distribution toward
$L(X^{\mathrm{exc}}
)$. In
Section \ref{sec4} we start by introducing the continuous-time height
process $H^{\mathrm{exc}}$ associated to $X^{\mathrm{exc}}$ and we
then show that $L(X^{\mathrm{exc}})$ can be
coded by $H^{\mathrm{exc}}$. In Section \ref{sec5} we use the
time-reversibility of $H^{\mathrm{exc}}$
to calculate the Hausdorff dimension of the stable lamination.

Throughout this work, the notation $\overline A$ stands for the closure
of a subset $A$ of the plane.

\begin{figure}

\includegraphics{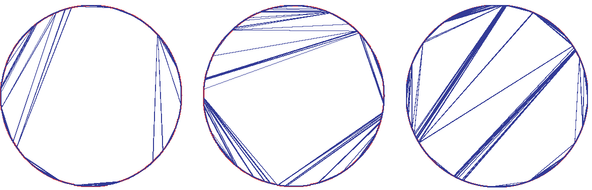}

\caption{Random dissections of $P_{27183}$ for $\theta
=1.1$, of $P_{11655}$
for $\theta=1.5$ and of $P_{20999}$ for $\theta=1.9$.}\label{figdiss}
\end{figure}

\section{The discrete setting: Dissections and trees}\label{sec1}

\subsection{Boltzmann dissections}\label{secb}

\begin{defn}A \emph{dissection} of a polygon is the union of the sides
of the
polygon and of a collection of diagonals that may intersect only at
their endpoints. A \emph{face} $f$ of a dissection $\omega$ of a
polygon $P$ is a connected component of the complement of $\omega$
inside $P$; its degree, denoted by $\deg(f)$, is the number of sides
surrounding $f$.
See Figure \ref{figdiss} for an example.
\end{defn}

Let $(\mu_i)_{i \geq2}$ be a sequence of nonnegative real numbers.
For every integer $n\geq3$,
let $P_n$ be the regular polygon of the plane whose vertices are the
$n$th roots of unity. For every $n \geq2$, let $\L_n$ be the set of
all dissections of
$P_{n+1}$. Note that $\L_n$ is a finite set. Let $\L= \bigcup_{n \geq
2} \L_n$ be the set of all dissections. A~weight $\pi(\omega)$ is
associated to each dissection $\omega\in\L_n$ by setting
\[
\pi(\omega)=\prod_{f\ \mathrm{ face\ of }\ \omega} \mu_{\deg(f)-1}.
\]
We define a probability measure on $\L_n$
by normalizing these weights. More precisely, we set
%
\begin{equation}
\label{eqzn}Z_n=\sum_{w \in\L_n} \pi(w),
\end{equation}
and for every $n \geq2$ such
that $Z_n > 0$,
\[
\P^ \mu_n(\omega)=\frac{1}{Z_n} \pi(\omega)
\]
for $\omega\in\L_n$.

We are interested in the asymptotic behavior of random dissections
sampled according to $\P^\mu_n$. Let $\Db$ be the closed unit disk of
the complex plane and let~$\C$ be the set of all compact subsets of
$\Db$.
We equip $\C$ with the Hausdorff distance $d_H$, so that $(\C,d_H)$
is a compact metric space. In the following, we will always view a
dissection as an element of this metric space.\vadjust{\goodbreak}

We are interested in the following question. For every
$n \geq2$ such that $Z_n>0$, let~$\mathfrak l_n$ be a random dissection
distributing according
to $\P^\mu_n$. Does there exist a limiting random compact set
$\mathfrak l$ such
that $\mathfrak l_n$ converges in distribution toward $\mathfrak l$?

We shall answer this question for some specific families of sequences
$(\mu_i)_{i \geq2}$ defined as follows. Let $ \theta\in(1, 2]$. We
say that a
sequence of nonnegative real numbers $(\mu_j)_{j \geq2}$ satisfies
the condition $(H_\theta)$ if:
%
\begin{longlist}[$-$]
\item[$-$] $\mu$ is
critical, meaning that $\sum_{i=2}^{\infty} i \mu_i =1$. Note that
this condition implies $\sum_{i=2}^{\infty} \mu_i<1$.
\item[$-$] Set $\mu_1=0$ and $\mu_0=1-\sum_{i=2}^{\infty}
\mu_i$. Then $(\mu_j)_{j \geq0}$ is a
probability measure in the domain of attraction of a stable law of
index $\theta$.
\end{longlist}

Recall that the second condition is equivalent to saying that if
$X$ is a random variable such that $\P[X=j]=\mu_j$ for $j \geq0$,
then either $X$ has finite variance or $\P[X \geq j]= j^{-\theta} L(j)$,
where $L$ is a function such that $\lim_{x \rightarrow\infty}
L(tx)/L(x)=1$ for all $t>0$ (such a function is called slowly
varying at infinity). We refer to \cite{Bingham} or \cite{Durrett}, Chapter
3.7, for details.

\subsection{Random dissections and Galton--Watson trees}\label{sec1.2}

In this subsection we explain how to associate a dual object to a
dissection. This dual object is a finite rooted ordered tree. The
study of large random dissections will then boil down to the study
of large Galton--Watson trees, which is a more familiar realm.

\begin{defn}Let $\N=\{0,1,\ldots\}$ be the set of all nonnegative
integers, $\N^*=\{1,2,\ldots\}$, and let $U$ be the set of labels
\[
U=\bigcup_{n=0}^{\infty} \bigl(\N^*
\bigr)^n,
\]
where by convention $(\N^*)^0=\{\varnothing\}$. An element of $U$ is
a sequence $u=u_1 \cdots u_m$ of positive integers, and we set
$|u|=m$, which represents the ``generation'' of $u$. If $u=u_1
\cdots u_m$ and $v=v_1 \cdots v_n$ belong to $U$, we write $uv=u_1
\cdots u_m v_1 \cdots v_n$ for the concatenation of $u$ and $v$. In
particular, note that $u \varnothing=\varnothing u = u$. Finally, a
\emph{rooted ordered tree} $\tau$ is a finite subset of $U$ such
that:
\begin{longlist}[(3)]
\item[(1)] $\varnothing\in\tau$;
\item[(2)] if $v \in\tau$ and $v=uj$ for some $j \in\N^*$, then $u
\in\tau$;
\item[(3)] for every $u \in\tau$, there exists an integer $k_u(\tau)
\geq0$ such that, for every $j \in\N^*$, $uj \in\tau$ if and only
if $1 \leq j \leq k_u(\tau)$.
\end{longlist}
In the following, by \emph{tree} we will always mean rooted ordered
tree. We denote the set of all trees by $\T$. We will often view
each vertex of a tree $\tau$ as an individual of a population whose
$\tau$ is the genealogical tree. The total progeny of $\tau$,
Card$(\tau)$, will be denoted by $\zeta(\tau)$. A leaf of a tree
$\tau$ is a vertex $u \in\tau$ such that $k_u(\tau)=0$. The total
number of leaves of $\tau$ will be denoted by $\lambda(\tau)$. If
$\tau$ is a tree and $u \in\tau$, we define the shift of $\tau$ at
$u$ by $T_u \tau=\{v \in U; uv \in\tau\}$, which is itself a
tree.
\end{defn}

\begin{figure}

\includegraphics{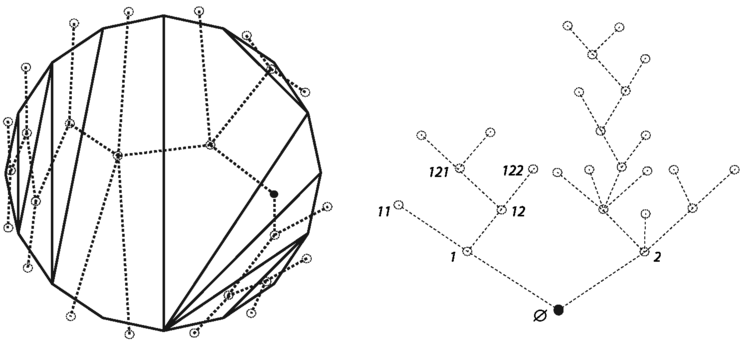}

\caption{The dual tree of a dissection, rooted at the bold
vertex.}\label{figdual}
\end{figure}

Given a dissection $\omega\in\L_n$, we construct a (rooted ordered)
tree $\phi(\omega)$ as follows: consider the dual graph
of $\omega$, obtained by placing a vertex inside each face of
$\omega$ and outside each side of the polygon $P_{n+1}$ and by
joining two vertices if the corresponding faces share a common edge,
thus giving a connected graph without cycles. Then remove the dual
edge intersecting the side $ [1,e^{{2 \textnormal{i} \pi
}/{(n+1)}} ]$
of $P_n$. Finally, root the graph at the dual vertex corresponding to
the face adjacent to the side $ [1,e^{{2i
\pi}/{(n+1)}} ]$ (see Figure \ref{figdual}). The planar structure
now allows us to associate a tree $ \phi( \omega)$ to this graph, in a
way that should be obvious from Figure \ref{figdual}. Note that $
k_u( \phi( \omega)) \neq1$ for every $u \in\phi( \omega)$.

For every integer $n \geq2$, let $ \T_ {(n)}$ stand for the set of all
trees $ \tau\in\T$ with exactly $n$ leaves and such that $k_u(\tau)
\neq1$ for every $u \in\tau$. The preceding construction provides a
bijection $ \phi$ from $ \L_n$ onto $ \T_ {(n)}$. Furthermore, if $
\tau= \phi( \omega)$ for $ \omega\in\L_n$, there is a one-to-one
correspondence between internal vertices of $ \tau$ and faces of $
\omega$, such that if $u$ is an internal vertex of $ \tau$ and $f$ is
the associated face of $ \omega$, we have $ \deg f = k_u( \tau)+1$. The
latter property should be clear from our construction.

%

\begin{defn}\label{defGW}Let $\rho$ be a probability measure on $\N$
with mean less than or equal to $1$ and such that $\rho(1)<1$. The law
of the
Galton--Watson tree with offspring distribution $\rho$ is the unique
probability measure $\P_\rho$ on $\T$ such that:\vadjust{\goodbreak}
\begin{longlist}[(1)]
\item[(1)] $\P_\rho[k_\varnothing=j]=\rho(j)$ for $j \geq0$;
\item[(2)] for every $j \geq1$ with $\rho(j)>0$, the shifted trees
$T_1 \tau, \ldots, T_j \tau$ are independent under the conditional
probability $\P_\rho[ \cdot| k_\varnothing=j]$ and their
conditional distribution is $\P_\rho$.
\end{longlist}
\end{defn}

A random tree with distribution $\P_\rho$ will sometimes be called a
$GW_\rho$ tree.


\begin{prop}\label{propGW}Let $(\mu_j)_{j \geq2}$ be a sequence of
nonnegative real numbers such that $\sum_{j=2}^{\infty} j\mu_j=1$. Put
$ \mu_1=0$ and $\mu_0=1-\sum_{j=2}^{\infty} \mu_j$ so that $ \mu=
( \mu_j)_ { j \geq0}$ defines a probability measure on $ \N$, which satisfies
the assumptions of Definition \ref{defGW}. Let $ n \geq2$ and let
$Z_n$ be defined as in \eqref{eqzn}. Then $Z_n > 0$ if, and only if, $
\mathbb{P}_{\mu}[\lambda(\tau)=n]>0$. Assume that this condition
holds. Then if
$\mathfrak l_n$ is
a random dissection distributed according to $\P^\mu_n$, the
tree $\phi(\mathfrak l_n)$ is distributed
according to $\mathbb{P}_{\mu}[\cdot| \lambda(\tau)=n]$.
\end{prop}

\begin{pf} Let $ \tau\in\T_ {(n)}$ and $ \omega= \phi^ {-1}( \tau
)$. Then
%
\begin{equation}
\label{eqGW} \Pmu(\tau)= \prod_ {u \in\tau}
\mu_ {k_u ( \tau)} = \mu_0^n \prod
_ {f\ \mathrm{ face\ of }\ \omega} \mu_ { \deg(f)-1}= \mu_0^n
\pi( \omega).
\end{equation}
The first equality is a well-known property of Galton--Watson trees
(see, e.g., Proposition 1.4 in \cite{RandomTrees}). The second one
follows from the observations preceding Definition \ref{defGW}, and
the last one is the definition of $ \pi( \omega)$. From \eqref{eqGW},
we now get that $\Pmu( \T_ {(n)})= \mu_0^n Z_n$, and then (if these
quantities are positive) that $\Pmu( \tau\mid\T_ {(n)})= \P^ \mu_n(
\omega)$, giving the last assertion of the proposition.
\end{pf}

\begin{rem} The preceding proposition will be a major ingredient of our
study. We will derive information about the random dissection $
\mathfrak l_n$
(when $n \rightarrow\infty$) from asymptotic results for the random
trees $ \phi( \mathfrak l_n)$. To this end, we will assume that $( \mu_j)_ {j
\geq2}$ satisfies condition $(H_ \theta)$ for some $ \theta\in
(1,2]$, which will allow us to use the limit theorems of \cite{K} for
Galton--Watson trees conditioned to have a (fixed) large number of leaves.
\end{rem}

\subsection{Coding trees and dissections}\label{sec1.3}

\begin{figure}

\includegraphics{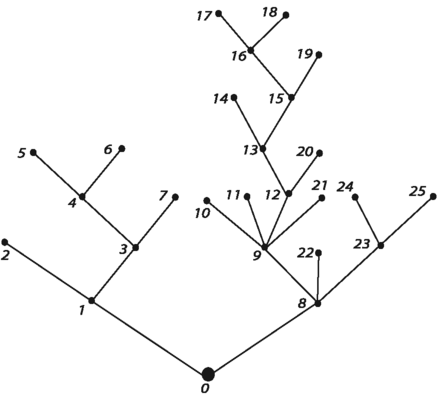}

\caption{The dual tree $\tau$ associated to the
dissection of Figure \protect\ref{figdual} with its vertices indexed in
lexicographical order. Here, $\zeta(\tau)=26$.}\label{figtree1}\vspace*{3pt}
\end{figure}

\begin{figure}[b]\vspace*{3pt}

\includegraphics{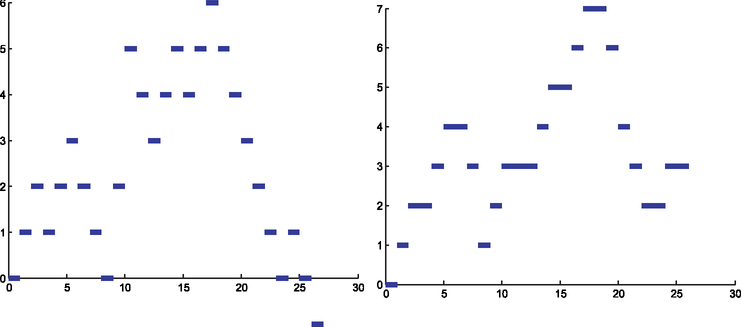}

\caption{The Lukasiewicz path $(W_{u}(\tau), 0 \leq u \leq\zeta
(\tau))$ and
the height function
$(H_{u}(\tau), 0 \leq u < \zeta(\tau)$ of $\tau$.}\label{figtree}
\end{figure}

In the previous subsection we have seen that certain random dissections
are coded by conditioned
Galton--Watson trees. We now
explain how trees themselves can be coded by two functions, called,
respectively, the
Lukasiewicz path and the height function (see Figures \ref{figtree1} and \ref{figtree}
for an example). These codings are crucial in the understanding of
large Galton--Watson trees and thus of large random dissections.

We write $u<v$ for the lexicographical order on the set $U$ (e.g.,
$\varnothing<1<21<22$). In the following, we will denote the children
of a tree $ \tau$ listed in lexicographical order by
$\varnothing=u(0)<u(1)<\cdots<u(\zeta(\tau)-1)$.

%

%

\begin{defn} Let $ \tau\in\T$. The height process
$H(\tau)=(H_n(\tau), 0 \leq n < \zeta(\tau))$ is defined, for $0
\leq n < \zeta(\tau)$, by $H_n(\tau)=|u(n)|$.
The Lukasiewicz path $W(\tau)=(W_n(\tau), 0 \leq n \leq\zeta(\tau))$
is defined by
$W_0(\tau)=0$ and $W_{n+1}(\tau)=W_{n}(\tau)+k_{u(n)}(\tau)-1$ for $0
\leq n \leq\zeta(\tau)-1$.
\end{defn}
 It is easy to see that $W_n ( \tau) \geq0$ for $0 \leq n <
\zeta(\tau)$ but $W_\zeta(\tau)=-1$ (see, e.g.,~\cite{RandomTrees}).

Consider a dissection $\omega$, its dual tree $\tau=\phi(\omega)$
and $W(\tau)$, the associated Lukasiewicz path. We now explain how to
reconstruct $ \omega$ from $W(\tau)$. As a
first step, recall that an internal vertex $u$ of $\tau$ is associated to
a face $f$ of $\omega$, and that the chords bounding $f$ are in
bijection with the dual edges linking $u$ to its children and to its
parent. The following proposition explains how to find all the children
of a given vertex of $\tau$ using only $W$ or $H$, and will be useful
to construct the edges linking the vertex $u \in\tau$ to its children.

\begin{prop}\label{propdiscretesons}
Let $ \tau\in\T,$ and let $u(0), \ldots, u( \zeta(\tau)-1)$ be as
above the
vertices of $ \tau$ listed in lexicographical order. Fix $ n \in\{
0,1, \ldots, \zeta(\tau)-1\}$ such that $k_ {u(n)} ( \tau)>0$ and
set $k=k_ {u(n)}( \tau)$.
\begin{longlist}[(ii)]
\item[(i)] Let $s_1, \ldots, s_k \in\{0,1, \ldots, \zeta(\tau)-1\}
$ be defined
by setting $s_i= \inf\{l \geq n+1 ; W_l ( \tau)=W_ {n+1} ( \tau
)-(i-1)\}$ for $1 \leq i \leq k$ (in particular, $s_1=n+1$). Then
$u(s_1),u(s_2), \ldots, u(s_k)$ are the children of $u(n)$ listed in
lexicographical
order.
\item[(ii)] We have
$H_ {s_1} ( \tau)=H_ {s_2}( \tau)=\cdots=H_ {s_k}( \tau)=H_n( \tau)+1$.
Furthermore, for $1 \leq i \leq k-1$,
\[
H_j( \tau)> H_ {s_i}( \tau)=H_ {s_ {i+1}} ( \tau)\qquad
\forall j \in( s_r,s_{r+1} ) \cap\N.
\]
\end{longlist}
\end{prop}

\begin{pf}We leave this as an exercise (or see the proof of Proposition
1.2 in~\cite{RandomTrees}) and encourage the reader to visualize what
this means on Figure \ref{figtree}.
\end{pf}

In a second step, we explain how to reconstruct the dissection from
the Lukasiewicz path of its dual tree.

\begin{prop}\label{deflien}Let $\zeta\geq2$ be an integer and let
$Z=(Z_n, 0 \leq n \leq\zeta)$ be a sequence of integers such that
$Z_0=0$, $Z_\zeta=-1$, $Z_k \geq0$ for $0 \leq k < \zeta$ and
$Z_{i+1}-Z_i \in\{ -1,1,2,3,\ldots\}$ for $0 \leq i < \zeta$. For
$0 \leq i <\zeta$, set $X_i=Z_{i+1}-Z_i$ and, for $ 1 \leq i \leq
\zeta$,
\[
\Lambda(i)=\card\{ 0 \leq j < i; X_j=-1\}.
\]
For every integer $i\in\{0,1,\ldots,\zeta(\tau)-1\}$ such that
$X_i \geq1$, set $k_i=X_i+1$ and let
$s^i_1,\ldots,s^i_{k_i+2}$ be defined by $s^i_1=s^i_{k_i+2}=i+1$ and
$s^i_{m+1}= \inf\{ l\geq
i+1 ; Z_l=Z_{i+1}-m \}$ for $1 \leq m \leq k_i$. Then the set $D(Z)$ defined
by
%
\begin{equation}
\label{eqDZ}D(Z)=\bigcup_{ i ; X_i \geq1} \bigcup
_{j=1}^{k_i+1} \biggl[\exp \biggl(-2 \textnormal{i} \pi
\frac{\Lambda(s_j^i)}{\Lambda(\zeta)+1 } \biggr), \exp \biggl(-2 \textnormal{i} \pi \frac{\Lambda(s_{j+1}^i)}{\Lambda(\zeta)+1 }
\biggr) \biggr]
\end{equation}
is a dissection of the polygon $P_{\Lambda(\zeta)+1}$ called the
dissection coded by $Z$.
\end{prop}

Note that if $\tau$ is a tree (different from the trivial tree $\{
\varnothing\}$), if $u(0), \ldots,$ $u( \zeta(\tau)-1)$ are its
vertices listed
in lexicographical order and $Z=W(\tau)$, then $\Lambda(i)$\vadjust{\goodbreak}
is the number of leaves among $ u(0), u(1), \ldots u(i-1)$ [in
particular, $\Lambda( \zeta)$ is the number of leaves of $ \tau$],
$k_i$ is the number of children of $u(i)$, and $s^i_m$ is the index of
the $m$th child of $u(i)$ for $1 \leq m \leq k_i$.

\begin{pf} First notice that, for all pairs $(i,j)$ occurring in the
union of \eqref{eqDZ}, we have $ \Lambda( s^i_j) \neq\Lambda(s^i_ {j+1})$.
We then
check that all edges of the polygon $P_{\Lambda(\zeta)+1}$ appear in
the right-hand side
of (\ref{eqDZ}). To this end, fix $\ell\in\{0,1,\ldots, \Lambda
(\zeta
)-1\}$. Then there is a unique integer
$k\in\{1,2,\ldots,\zeta-1\}$ such that $X_k=-1$ and $\Lambda
(k)=\ell$. Set
\[
i=\sup\bigl\{j\in\{0,1,\ldots,k-1\}\dvtx Z_j\leq Z_k
\bigr\}
\]
and $m=Z_{i+1}-Z_k+1$. Notice that $1\leq m\leq k_i$ since $Z_k\geq
Z_i$ by construction.
It is now a simple matter to verify that $s^i_m=k$ and $s^i_{m+1}=k+1$.
Recalling that $\Lambda(k)=\ell$
and $\Lambda(k+1)=\ell+1$, we get that the line segment
\[
\biggl[\exp \biggl(-2 \textnormal{i} \pi \frac{\ell}{\Lambda(\zeta)+1 } \biggr), \exp
\biggl(-2 \textnormal {i} \pi \frac{\ell+1}{\Lambda(\zeta)+1 } \biggr) \biggr]
\]
appears in the right-hand side
of (\ref{eqDZ}). We therefore get that $D(Z)$ contains all edges of
$P_{\Lambda(\zeta)+1}$
with the possible exception of the edge $[1,\exp(-2 \textnormal{i}
\pi
\frac{\Lambda(\zeta)}{\Lambda(\zeta)+1 })]$. However, the latter edge
also appears in the union
of (\ref{eqDZ}), taking $i=0$ and $j=k_0+1$ and noting that
$s^0_{k_0+1}= \zeta$ and $s^0_{k_0+2}=1$.

Next suppose that $0 \leq i < \z, 0 \leq i' < \zeta$ are such that $k_i
\geq1, k_ {i'} \geq1$. Let $ j \in\{1, \ldots, k_i+1\}$, $ j' \in
\{1, \ldots, k_ {i'}+1\}$. If $(i,j) \neq(i',j')$, one easily checks
that either the intervals $(s^i_j, s^i_ {j+1})$ are disjoint or one of
them is contained in either one. It follows that the chords
corresponding, respectively, to $(i,j)$ and to $(i',j')$ in the union of
\eqref{eqDZ} are noncrossing. Hence, $D(Z)$ is a dissection.
\end{pf}

\begin{lem}\label{lemecriture}For every dissection $\omega\in\L$, we
have $D(W(\phi(\omega)))=\omega$.
In other words, a dissection is equal to the dissection coded by the
Lukasiewicz path of its dual tree.
\end{lem}

\begin{pf}This is a consequence of our construction. Suppose that $
\omega\in\L_n$,
for some $n\geq2$, and set $\tau=\phi(\omega)$. Fix a face $f$ of
$\omega$ and the corresponding
dual vertex $u(i) \in\phi(\omega)$ (recall that the faces of $f$ are
in one-to-one correspondence with the internal vertices of $\tau$).
Denote the Lukasiewicz path of $\tau$ by
$Z=W(\tau)$. First notice that the
degree of $f$ is equal to $1+k_{u(i)}=Z_{i+1}-Z_i+2$, where
$k_{u(i)}$ is the number of children of $u(i)$. To simplify
notation, set $k_i=k_{u(i)}$. Let $s^i_1,\ldots,s^i_{k_i+2}$ be defined
as in Proposition \ref{deflien}. By Proposition
\ref{propdiscretesons}, $u(s^i_1),u(s^i_2), \ldots, u(s^i_ {k_i})$
are the
children of $u(i)$.

As in Proposition \ref{deflien}, we set, for every $1 \leq i \leq
\zeta$, $\Lambda(i)=\card\{ 0 \leq j < i; Z_{j+1}-Z_j=-1\}$,
which represents the number of leaves among the first $i$ vertices
of $\tau$. Note that $ \Lambda( \zeta(\tau))=n$. Then, assuming
that $k_i
\geq2$:
\begin{longlist}[$-$]
\item[$-$] For every $1 \leq j \leq k_i$ the chord of $\omega$ which
intersects the dual edge linking
$u(i)$ to its $j$th child is
\[
\biggl[\exp \biggl(-2 \textnormal{i} \pi \frac{\Lambda(s_j^i)}{n+1 } \biggr), \exp
\biggl(-2 \textnormal{i} \pi \frac{\Lambda(s_{j+1}^i)}{n+1 } \biggr) \biggr].
\]
\item[$-$] The chord of $\omega$ intersecting the dual edge linking
$u(i)$ to its parent is
\[
\biggl[\exp \biggl(-2 \textnormal{i} \pi \frac{\Lambda(s_{ {k_i}+1}^i)}{n+1 } \biggr), \exp
\biggl(-2 \textnormal {i} \pi \frac{\Lambda(s_{1}^i)}{n+1 } \biggr) \biggr].
\]
\end{longlist}
Indeed, a look at Figure \ref{figdual} should convince the reader
that the vertices
\[
\exp \biggl(-2 \textnormal{i} \pi \frac{\Lambda(s_j^i)}{n+1 } \biggr),\qquad 1 \leq j \leq
k_i+1
\]
are exactly the vertices belonging to the boundary of the face
associated with $u(i)$ listed in clockwise order. Consequently, the
preceding chords are exactly the ones that bound this face. Since this
holds for every face $f$ of $ \omega$, the conclusion follows.
\end{pf}


\section{The continuous setting: Construction of the stable lamination}\label{sec2}

In this section we present the continuous background by
first introducing the normalized excursion $X^{\mathrm{exc}}$ of the
$\theta$-stable
L\'{e}vy process. This process is important for our purposes because
$X^{\mathrm{exc}}$
will appear as the limit in the Skorokhod sense of the rescaled
Lukasiewicz paths of large $GW_\mu$ trees coding discrete
dissections. We then use $X^{\mathrm{exc}}$ to construct a random
compact subset
of the closed unit disk, which will be our candidate for the limit in
distribution of the random dissections we are considering. Two cases
will be distinguished: the case $\theta=2$, where $X^{\mathrm{exc}}$ is
continuous, and the case $\theta\in(1,2)$, where the set of
discontinuities of $X^{\mathrm{exc}}$ is dense.

\subsection{The normalized excursion of the L\'{e}vy process}\label{sec2.1}

We follow the presentation of \cite{Duquesne} and refer to \cite
{Bertoin} for the proof of the results recalled in this subsection. The
underlying
probability space will be denoted by $(\Omega, \mathcal{F}, \P)$.
Let $X$ be a process with paths in $\D(\R_+, \R)$, the space of
right-continuous with left limits (c\`{a}dl\`{a}g) real-valued functions,
endowed with the Skorokhod topology. We refer the reader to
\cite{Bill}, Chapter 3 and \cite{Shir}, Chapter VI, for background
concerning the Skorokhod topology. We denote by $(\F_t)_{t \geq0}$
the canonical filtration of $X$ augmented with the
$\P$-negligible sets. We assume that $X$ is a strictly stable spectrally
positive L\'{e}vy process of index $\theta$ normalized so that for every
$\lambda>0$,
\[
\E\bigl[\exp(-\lambda X_t)\bigr]=\exp\bigl(t \lambda^\theta
\bigr).
\]
In the following, by the \emph{$ \theta$-stable L\'{e}vy process} we
will always mean
such a L\'{e}vy process. In
particular, for $\theta=2$ the process $X$ is $\sqrt{2}$ times the
standard Brownian motion on the line. Recall that $X$ enjoys the
following scaling property: For every $c>0$, the process $(c^{-1/\theta}
X_{ct}, t \geq0$) has the same law as $X$. Also recall that when $1
< \theta< 2$, the L\'{e}vy measure $\pi$ of $X$ is
\[
\pi(dr)=\frac{\theta(\theta-1)}{\Gamma(2-\theta)} r^{-\theta-1} 1_
{(0,\infty)}\,dr.
\]
For $s> 0$, we set $\Delta X_s= X_s-X_{s-}$. The following notation
will be useful: for $0 \leq s<t$,
\[
I^s_t = \inf_{[s,t]} X,\qquad I_t=
\inf_{ [0,t ]} X,\qquad S_t=\sup_{ [0,t ]} X.
\]
Notice that the process $I$ is continuous since $X$ has no negative jumps.

We have $X_0=0$ and $I_t < 0 < S_t$ for every $t >0$ almost surely
[meaning that the point $0$ is regular both for $(0, \infty)$ and for
$(- \infty,0)$ with respect to $X$]. The process $X-I$ is a strong
Markov process and $0$ is regular for
itself with respect to $X-I$. We may and will choose $-I$ as the
local time of $X-I$ at level $0$. Let $(g_i,d_i), i \in\mathcal{I}$
be the excursion intervals of $X-I$ away from $0$. For every $i \in
\mathcal{I}$ and $ s \geq0$, set $\omega_s^i= X_{(g_i+s) \wedge
d_i}-X_{g_i}$.
We view $ \omega^i$ as an element of the excursion space $ \mathcal
{E}$, which is defined by
\[
\mathcal{E}= \bigl\{ \omega\in\D(\R_+, \R_+); \omega(0)=0 \mbox{ and } \z (
\omega):= \sup\bigl\{s>0; \omega(s)>0 \bigr\} \in(0, \infty) \bigr\}.
\]
If $\omega\in\mathcal{E}$, we call $ \zeta( \omega)$ the lifetime of
the excursion $ \omega$. From It\^{o}'s excursion theory, the point measure
\[
\mathcal{N}( dt\, d \omega)= \sum_{i \in\mathcal{I}}
\delta_{(-I_{g_i},\omega^i)}
\]
is a Poisson measure with intensity $dt N(d\omega)$, where $N(d
\omega)$ is a $\sigma$-finite measure on the set $ \mathcal{E}$.

\begin{figure}

\includegraphics{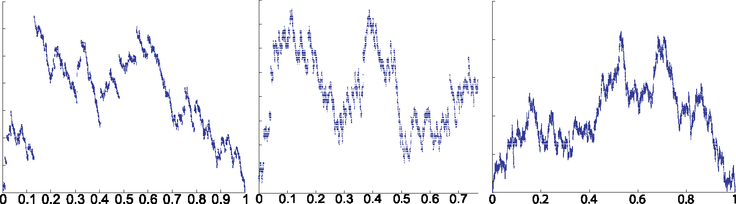}

\caption{Simulations of $X^{\mathrm{exc}}$ for, respectively, $\theta
=1.1,1.5,1.9$.}\label{figexc}
\end{figure}

Let us define the normalized excursion of the $\theta$-stable L\'{e}vy process.
Define, for every $\lambda>0$, the re-scaling operator $S^{(\lambda)}$
on the set of excursions by $S^{(\lambda)}(\omega)= ( \lambda^{1/\theta} \omega(s/\lambda), s \geq0  )$.
The scaling property of X shows that the image of $N(\cdot|
\zeta>t)$ under $S^{(1/\zeta)}$ does not depend on $t>0$. This
common law, which is supported on the c\`{a}dl\`{a}g paths with unit
lifetime, is called the law of the normalized excursion of $X$ and
denoted by $\P^\exc$. Informally, $\P^\exc$ is the law of an
excursion under the It\^{o} measure conditioned to have unit lifetime.
In the following, $(X^{\mathrm{exc}}_t; 0 \leq t \leq1)$ will stand
for a
process defined on $(\Omega,\F,\P)$ with paths in $\D([0,1],\R_+)$
and whose distribution under $\P$ is $\P^\exc$ (see Figure \ref{figexc} for a simulation).
Note that $
X^{\mathrm{exc}}_0=X^{\mathrm{exc}}_1=0$.

As for the Brownian excursion, the normalized excursion can be
constructed directly from the L\'{e}vy process $X$. We state Chaumont's
result \cite{Chaumont} without proof. Let
$(\underline{g}_1,\underline{d}_1)$ be the excursion interval of
$X-I$ straddling $1$. More precisely, $\underline{g}_1=\sup\{s \leq1;
X_s=I_s \}$ and $\underline{d}_1=\inf\{s > 1; X_s=I_s\}.$
Let $\zeta_1=\underline{d}_1-\underline{g}_1$ be the length of this
excursion.

\begin{prop}[(Chaumont)]\label{proplaw}Set $X^\ast_s=\zeta_1^{-1/\theta}
(X_{\underline{g}_1+\zeta_1 s} - X_{\underline{g}_1})$ for every $ s
\in[0,1]$.
Then $X^ \ast$ is distributed according to $\P^ \exc$.
\end{prop}

\subsection{\texorpdfstring{The $\theta$-stable lamination of the disk}
{The theta-stable lamination of the disk}}\label{sec2.2}

The open unit disk of the complex plane $\mathbb{C}$ is denoted by
$\D=\{z \in\mathbb{C}; |z|<1\}$ and $\mathbb{S}_1$ is the unit
circle. If $x,y$ are
distinct points of $\mathbb{S}_1$, we recall that $[x,y]$ stands for
the line segment between $x$ and $y$. By convention, $[x,x]$ is equal to
the singleton $\{x\}$.

\begin{defn}A geodesic lamination $L$ of $\Db$ is a closed subset $L$
of $\Db$ which can be written
as the union of a collection of noncrossing chords. The lamination
$L$ is maximal if it is maximal for the inclusion relation among
geodesic laminations of $\Db$. In the sequel, by lamination we will
always mean geodesic lamination of $\Db$.
\end{defn}

\begin{rem}In hyperbolic geometry, geodesic laminations of the disk are
defined as closed subsets of the open
hyperbolic disk \cite{Bon}. As in \cite{CLG}, we prefer to see these
laminations as compact subsets of $\Db$ because this will allow us
to study the convergence of laminations in the sense of the
Hausdorff distance on compact subsets of $\Db$.
\end{rem}

It is not hard to check that the set of all geodesic laminations is
closed with respect to the Hausdorff distance.

\subsubsection{The Brownian triangulation}

\begin{defn}The Brownian excursion $\mathbh{e}$ is defined
as $X^{\mathrm{exc}}$ for $\theta=2$. For $u,v \in[0,1]$ we set $u
\stackrel{\mathbh{e}}{\thicksim}v$ if $\mathbh{e}_{u\wedge
v}=\mathbh{e}_{u\vee
v}=\mathop{\min}_{t\in[u\wedge v,u\vee v]}\mathbh{e}_t.$
\end{defn}

Note that, with our normalization of $X^{\mathrm{exc}}$, $\mathbh{e}
/ \sqrt{2}$
is the standard Brownian excursion. It is well known that the local
minima of
$\mathbh{e}$ are distinct almost surely. In the following, we always discard
the set of probability zero where this property fails.

\begin{figure}

\includegraphics{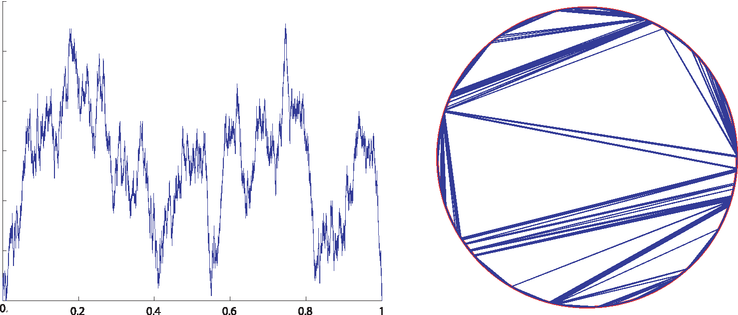}

\caption{A Brownian excursion $\me$ and the
associated triangulation $\mathbf{L}(\me)$.}\label{figtri}
\end{figure}

\begin{prop}[(Aldous \cite{Aldous}--Le Gall and Paulin \cite
{LGP})]\label
{propLe}Define
$\mathbf{L}(\mathbh{e})$ by
\[
\mathbf{L}(\mathbh{e})=\bigcup_{s \stackrel{\mathbh
{e}}{\thicksim
} t}
\bigl[e^{-2 \textnormal{i} \pi s},e^{-2 \textnormal{i} \pi
t}\bigr].
\]
Then $\mathbf{L}(\mathbh{e})$ is a maximal geodesic
lamination of $\Db$ [see Figure \ref{figtri} for a simulation of
$\mathbf{L}(\me)$].
\end{prop}

%

\begin{rem} Both the property that $\mathbf{L}(\mathbh{e})$ is a
lamination and its maximality property are related to the fact that
local minima of
$\mathbh{e}$ are distinct. The connected components
of $\Db\setminus\mathbf{L}(\mathbh{e})$ are open triangles whose
vertices belong to $\mathbb{S}_1$. For this reason we call
$\mathbf{L}(\mathbh{e})$ the Brownian triangulation. Notice also
that $\mathbb{S}_1 \subset\mathbf{L}(\mathbh{e})$.
\end{rem}

\subsubsection{\texorpdfstring{The $\theta$-stable lamination}
{The theta-stable lamination}}

Here, $\theta\in(1,2)$ so that the $\theta$-stable L\'{e}vy process
$X$ is not continuous.
In the beginning of this section we fix $Z \in\D([0,1], \R)$ such that
$Z_0=Z_1=0$, $\Delta Z_s \geq0$ for $s \in(0,1]$ and $Z_s >0$ for $s
\in(0,1)$. We then consider the case when $Z=X^{\mathrm{exc}}$ is the
normalized
excursion of the $\theta$-stable L\'{e}vy process
$X$.

\begin{defn}For $0 \leq s < t \leq1$, we set $ s \mathop{\simeq}^{Z}
t$ if and only if $t=
\inf\{ u > s; Z_u \leq Z_{s-}\}$ (where $Z_ {0-}=0$ by definition).
For $0 \leq t < s \leq1$, we
set $s \mathop{\simeq}^{Z} t$ if and only if $t \mathop{\simeq}^{Z}
s$. Finally,
we set $s \mathop{\simeq}^{Z} s$ for every $s \in[0,1]$.
\end{defn}

Note that
$\mathop{\simeq}^{Z}$ is not necessarily an equivalence relation. For
example, if $0 < r < s < t <1$ are such that $ \Delta Z_r=0$,
$Z_r=Z_s=Z_t$ and $Z_u> Z_r$ for $u \in(r,s) \cup(s,t)$, then $r
\mathop{\simeq}^{Z} s$ and $s \mathop{\simeq}^{Z} t$, but we do not
have $r \mathop{\simeq}^{Z} t$.


\begin{rem}\label{remeq}If $s \mathop{\simeq}^{Z} t$ and $s<t$, then
$Z_{s-}= Z_t$ and $Z_r > Z_{s-}$ for $r \in(s,t)$.
\end{rem}

\begin{prop}\label{proplamination}We say that $Z$ attains a local
minimum at $t \in
(0,1)$ if there exists $\eta>0$ such that $\inf_{ [t-\eta
,t+\eta ]} Z=Z_{t}$.
Suppose that $Z$ satisfies the following four assumptions:
\begin{longlist}[(H1)]
\item[(H1)] If $0 \leq s < t \leq1$, there exists at most one value
$r \in(s,t)$ such that $Z_r=\inf_{[s,t]} Z$ (we say
that local minima of $Z$ are distinct);
\item[(H2)] If $t \in(0,1)$ is such that $\Delta Z_t >0$, then $\inf_{ [t,t+\varepsilon ]} Z < Z_t$ for all
$0<\varepsilon\leq1-t$;
\item[(H3)] If $t \in(0,1)$ is such that $\Delta Z_t >0$, then $\inf_{ [t-\varepsilon,t ]} Z < Z_{t-}$ for all $ 0 < \varepsilon\leq t$;
\item[(H4)] Suppose that $Z$ attains a local minimum at $t \in(0,1)$
[in particular, $\Delta Z_t=0$ by \textup{(H3)}].
Let $s = \sup\{ r \in[0,t] ; Z_r < Z_{t}\}$. Then $\Delta Z_s
>0$ and $Z_{s-}<Z_{t}$. Note that then $Z_s>Z_t$ by \textup{(H2)}.
\end{longlist}
Then the set
\[
L(Z):= \bigcup_{s \mathop{\simeq}^Z t}\bigl[e^{-2 \textnormal{i} \pi
s}
,e^{-2 \textnormal{i} \pi t}\bigr]
\]
is a geodesic lamination of $\Db$, called the
lamination coded by the c\`{a}dl\`{a}g function~$Z$.
\end{prop}

Notice that $\mathbb{S}_1 \subset L(Z)$ since $s \mathop{\simeq}^{Z} s$
for every $s \in[0,1]$.

\begin{pf}It easily follows from Remark \ref{remeq} that the chords
appearing in the definition of $L(Z)$ are noncrossing. We have to prove
that $L(Z)$ is closed. To this end, it is enough to verify that the
relation $\mathop{\simeq}^{Z}$ is closed, in the sense that its graph
is a closed subset of $[0,1]^2$.
Consider two sequences
$(s_n),(t_n)$ of reals such that $0 \leq s_n<t_n \leq1$, $s_n
\mathop{\simeq}^Z t_n$ and the pairs $(s_n,t_n)$ converge to $(s,t)$.
We need to verify that $s \mathop{\simeq}^{Z} t$. Clearly, $s \leq
t$ and we can assume that $s<t$ since $\mathbb{S}_1 \subset L(Z)$.

The property $s_n
\mathop{\simeq}^Z t_n$ implies that $Z_r \geq Z_ {t_n}$ for every $r
\in(s_n,t_n)$. By passing to the limit $n \rightarrow\infty$, we get
$Z_r \geq Z_ {t-}$ for every $ r \in(s,t)$. If $ \Delta Z_t>0$, this
contradicts~(H3). So we can assume that $ \Delta Z_t=0$, implying that
the sequence $(Z_{t_n})$
converges to $Z_t$ as $n \rightarrow\infty$.

\emph{Case} 1. Assume that $\Delta Z_s
>0$ and thus $s>0$. By (H2) and right-continuity, we can find $\eta>0$ such
that $ \eta< (t-s)/2$ and
\[
\inf_{[s,s+\eta)} Z > \inf_{[s+\eta,(s+t)/2]} Z.
\]
It follows from (H3) that the infimum of $Z$ over a compact interval is
achieved at some point of this interval. Hence, we may take $r_0 \in
[s+\eta,(s+t)/2]$ such that
$Z_ {r_0}=\inf_{[s+\eta,(s+t)/2]} Z$. If $s<s_n$ for infinitely many $n$,
we can find infinitely many values of $n$ for which
$s<s_n<s+\eta\leq r_0<t_n$. For those values of $n$, $r_0 \in
(s_n,t_n)$ and
$Z_ {r_0} < Z_{s_n-}$, which contradicts Remark \ref{remeq}. We can
thus suppose that $s_n \leq s$ for every sufficiently large $n$.
Consequently, $(Z_{s_n-})$ converges
to $Z_{s-}$ as $n$ tends to infinity. Since $Z_{s_n-}=Z_{t_n}$ for
all $n$, it follows that $Z_t=Z_{s-}$.
Recall that $Z_r \geq Z_t$ for $r \in(s,t)$. We now demonstrate by
contradiction that, in fact, $Z_r > Z_t$ for
all $r \in(s, t)$. Suppose
that there exists $r_1 \in(s,t)$ such that $Z_ {r_1}=Z_t$. Notice that
$Z$ then attains a local minimum at $r_1$. Property (H3) ensures that
\[
s= \sup\bigl\{ u \in[0,r_1]; Z_u < Z_ {r_1}
\bigr\},
\]
and the fact that
$Z_{s-}=Z_{t}=Z_ {r_1}$ contradicts (H4). We conclude that $Z_ {r}
>Z_ {s-}$ for every $r \in(s,t)$. Therefore, $t=\inf\{ u > s ; Z_u
\leq Z_{s-} \}$. This implies
that $s \mathop{\simeq}^{Z} t$, as desired.

\emph{Case} 2. Assume that $\Delta Z_s = 0$. In this case, $(Z_{s_n})$
converges to $Z_s$ as $n$ tends to infinity. Since
$Z_{s_n-}=Z_{t_n}$ for all $n$, it follows that $Z_s=Z_t$. We also know
that $Z_r \geq Z_s$ for $r \in(s,t)$. If $s=0$, we necessarily have
$t=1$ and the fact that $Z$ is positive on $(0,1)$ implies $0 \mathop
{\simeq}^{Z}
1$. We thus suppose that $s>0$. Argue by contradiction and
suppose that there exists $r_1 \in(s,t)$ such that $Z_ {r_1}=Z_t$.
Then $r_1$ is a local minimum of $Z$. If
$\inf_{ [s-\varepsilon,s ]} Z < Z_s$ for every $ \varepsilon\in
(0,s]$, then
$s=\sup\{ u
\in[0,r_1] ;
Z_u < Z_ {r_1}\}$. By (H4), $s$ must be a jump time of $Z$, which is a
contradiction. If $\inf_{ [s-\varepsilon,s ]} Z\geq Z_s$ for
some $\varepsilon
\in(0,s]$, this means
that $s$ is a local minimum of $Z$. Since $Z_s=Z_ {r_1}$, this
contradicts (H1). We conclude that $Z_r
>Z_t$ for $r \in(s,t)$. This implies that $s \mathop{\simeq}^{Z}
t$.
\end{pf}

%

Let (H0) be the property: $\{ s \in[0,1]; \Delta Z_s \neq0
\}$ is dense in $[0,1]$.

\begin{prop}\label{propHX} Let $1 < \theta< 2$. With probability one,
the normalized excursion $X^{\mathrm{exc}}$ of the $\theta$-stable L\'
{e}vy process
satisfies the assumptions \textup{(H0), (H1), (H2), (H3)} and \textup{(H4)}.
\end{prop}

\begin{pf}
It is sufficient to prove that properties analogous to (H0)--(H4)
hold for the L\'evy process $X$. The case of (H0) is
clear. (H1) and (H2) are consequences of the (strong) Markov property
of $X$ and the fact that $0$ is regular for $(- \infty,0)$ with respect
to $X$.

For the remaining properties, we will use the time-reversal property
of $X$, which states that if $t>0$ and $\widehat{X}^{(t)}$ is the
process defined by
$\widehat{X}^{(t)}_s=X_t - X_{(t-s)-}$ for $0 \leq s <t$ and $\widehat
{X}^{(t)}_t=X_t$,
then the two processes $(X_s,0 \leq s \leq t)$ and $(\widehat{X}^{(t)}_s,0
\leq s \leq t)$ have the same law. For (H3), the time-reversal
property of $X$ and the regularity of~$0$ for $(0,\infty)$ shows
that a.s. for every jump time $s$ of $X$ and every $v \in[0,s)$,
\[
\inf_{r \in[v,s]}X_r < X_{s-}.
\]

We finally prove the analog of (H4) for $X$. By the time-reversal
property of $X$,
it is sufficient to prove that if $q>0$ is rational and $T= \inf\{t
\geq q; X_t \geq S_q\}$, then $X_T > S_q \geq X_ {T-}$ almost surely.
This follows from the Markov property at time~$q$ and the fact that
for any $a>0$, $X$ jumps a.s. across $a$ at its first passage time
above~$a$ (see \cite{Bertoin}, Proposition VIII.8 (ii)).
\end{pf}

In the following, we always discard the set of zero probability where
one of the
properties (H0)--(H4) does not hold.

\begin{defn}The $\theta$-stable lamination is defined as the geodesic
lamination $L(X^{\mathrm{exc}})$, where $X^{\mathrm{exc}}$ is the
normalized excursion of the
$\theta$-stable L\'{e}vy process.
\end{defn}

See Figure \ref{figdiss} for some examples. The following proposition
is immediate
from the definition of the relation $
\mathop{\simeq}^{X^\exc}$ and Remark \ref{remeq}.

\begin{prop} \label{propidentifications}Almost surely, for every
choice of $0\leq\a< \b\leq
1$ with $(\a,\b) \neq(0,1)$, we have $\a
\mathop{\simeq}^{X^\exc} \b$ if and only if one of the following two
mutually exclusive cases holds:
\begin{longlist}[(ii)]
\item[(i)] $\Delta X^{\mathrm{exc}}_\a>0$ and $\b= \inf\{ u \geq
\a;
X^{\mathrm{exc}}_u=X^{\mathrm{exc}}_{\a-}\}$;
\item[(ii)] $\Delta X^{\mathrm{exc}}_\a=0$, $X^{\mathrm{exc}}_\a
=X^{\mathrm{exc}}_\b$, and $X^{\mathrm{exc}}_r > X^{\mathrm
{exc}}_\a
$ for
every $r \in(\a,\b)$.
\end{longlist}
\end{prop}

\begin{defn}\label{defidentifications}Let $\mathcal{E}_{1}$ be the set
of all pairs $(\a,\b)$ where $0 \leq\a< \b\leq1$ satisfy
condition (i)
in Proposition \ref{propidentifications}.
\end{defn}

\begin{prop}\label{propapproxX}The following holds almost surely for
any pair $(s,t)$ such that $0 \leq s < t \leq1$
and $X^{\mathrm{exc}}_s=X^{\mathrm{exc}}_t$ and $X^{\mathrm{exc}}_r
> X^{\mathrm{exc}}_s$ for every $r\in(s,t)$. For
every $\varepsilon\in(0,(t-s)/2)$, there exist $s' \in
[s,s+\varepsilon]$ and $t' \in[t-\varepsilon,t]$ such that $\Delta
X^{\mathrm{exc}}_{s'}>0$ and $t'= \inf\{ u \geq s'; X^{\mathrm
{exc}}_u=X^{\mathrm{exc}}_{s'-}\}$, so that in
particular \mbox{$(s', t') \in\mathcal{E}_ {1}$}.
\end{prop}

\begin{pf} Let $ 0 \leq s<t \leq1$ be such that the assumptions in the
proposition hold. Take $\varepsilon< (t-s)/4$, then set
$m=\inf_{[s+\varepsilon,t-\varepsilon]}X^{\mathrm{exc}}$ and note that
$m>X^{\mathrm{exc}}_s$ as an easy
consequence of (H3).
By right-continuity, there exists $\varepsilon'$ with $0 <\varepsilon'<
\varepsilon$ such that $\sup_{[s,s+\varepsilon']}X^{\mathrm{exc}}< m$.
Let $w \in
(s,s+\varepsilon')$ be a jump time of $X^{\mathrm{exc}}$, so that, by (H2),
\[
\inf_{r \in[w, s +
\varepsilon']} X^{\mathrm{exc}}_r< X^{\mathrm{exc}}_w.
\]
We already noticed that the property (H3) implies that
the minimum of $X^{\mathrm{exc}}$ over a compact interval is achieved
at a point of
this interval. Hence,
there exists $u \in[w, s + \varepsilon]$ such that $X^{\mathrm{exc}}_u=
\inf_{[w, s+
\varepsilon]} X^{\mathrm{exc}}$. Finally, let $s'=\sup\{ r \in[s,u];
X^{\mathrm{exc}}_r<X^{\mathrm{exc}}_{u}\}$. By (H4), we see that $s'$
is a jump time. Set $t'=
\inf\{ u > s'; X^{\mathrm{exc}}_u=X^{\mathrm{exc}}_{s'-}\}$. By
construction, $s \leq
s' \leq w \leq u \leq s+\varepsilon< t-\varepsilon\leq t' \leq t$ and the
desired result follows.
\end{pf}

\begin{prop}\label{propfermL}We have a.s.
\[
L\bigl(X^{\mathrm{exc}}\bigr)= \overline{\bigcup_{(s,t) \in\mathcal{E}_{1}}
\bigl[e^{-2
\textnormal
{i} \pi
s},e^{-2 \textnormal{i} \pi t}\bigr]}.
\]
\end{prop}

\begin{pf}Denote the compact subset of $\Db$ in the right-hand side by
$K$. The fact that $L(X^{\mathrm{exc}})$ is closed implies that $K
\subset L(X^{\mathrm{exc}})$.
We have to show the reverse inclusion. To this end, let $0\leq u < v
\leq1$ such that $u
\mathop{\simeq}^{X^\exc} v$ but $(u,v)\notin\mathcal{E}_{1}$. Then
condition (ii) in Proposition
\ref{propidentifications} holds for $(\alpha,\beta)=(u,v)$, and it
follows from Proposition
\ref{propapproxX} that $(u,v)$ is the limit of a sequence of pairs
$(u_n,v_n)$ belonging to~$\mathcal{E}_{1}$.
Since $K$ is closed, we get that $[e^{-2 \textnormal{i} \pi u},e^{-2
\textnormal{i} \pi v}] \subset K$.
Finally, from the fact that $X^{\mathrm{exc}}$ satisfies properties
(H0) and (H2), it
is easy to verify that in any nontrivial open
subinterval of $[0,1]$ we can find a pair $(u,v)$ such that $(u,v)\in
\mathcal{E}_{1}$, and it follows that $\mathbb{S}_1 \subset K$. This completes
the proof.
\end{pf}


\section{Convergence to the stable lamination}\label{sec3}

In this section we show that the Boltzmann dissections of $ P_ {n+1}$
considered in Section \ref{secb} converge in distribution to the stable
laminations introduced in the previous section. To this end,
we use limit theorems for rescaled Lukasiewicz paths of critical
Galton--Watson trees conditioned on
their number of leaves, which we obtained in \cite{K}. We combine
these limit theorems with Proposition \ref{propGW} (which
states that the dual tree of a Boltzmann dissection is a
Galton--Watson tree conditioned on having a given number of leaves) to deduce
that the underlying tree structures of large dissections converge. As
before, we will deal separately with the case $\theta=2$ and the case
$\theta\in(1,2)$. Our
goal is to prove the following:

\begin{thmm}\label{thmcv1} Let $ ( \mu_j)_ { j \geq2}$ be a sequence
satisfying Assumption $(H_ \theta)$ for some $ \theta\in(1,2]$. For
every integer $n \geq2$ such that the definition of $\P^\mu_n$ makes
sense, let $\mathfrak l_n$ be a random dissection
distributed according to $\P^\mu_n$. Then
\[
\mathfrak l_n \mathop{\longrightarrow}^{(d)}_{n \rightarrow\infty}
\cases{ %
\mathbf{L}(\mathbh{e}), &\quad $\mbox{if }
\theta=2,$
\vspace*{2pt}\cr
 L\bigl(X^{\mathrm{exc}}\bigr),&\quad$\mbox{if }\theta\in(1,2),$}
\]
where the convergence holds in distribution for the Hausdorff distance
on the space of all
compact subsets of $\Db$.
\end{thmm}

%
\begin{rems}
(i) This theorem generalizes Aldous'
result \cite{Aldous,Aldous2}, stating that uniformly distributed
triangulations of
$P_n$ converge to $\mathbf{L}(\me)$ as $n \rightarrow\infty$. Indeed,
in our setting,
uniform triangulations of $P_n$ are obtained by taking
$\mu_0=1/2,\mu_2=1/2$ and $\mu_j = 0$ otherwise.\vspace*{-6pt}
\begin{longlist}[(iii)]
\item[(ii)]In \cite{CK}, it is shown that Theorem \ref{thmcv1} can be
used to study
uniformly distributed dissections. More precisely, if one sets $ \mu_0
= 2 - \sqrt{2}$ and $ \mu_i= ((2 - \sqrt{2})/2)^ {i-1}$ for every $i
\geq2$, then the Boltzmann probability measure $\P^\mu_n$
associated to $\mu$ is the uniform probability measure on
dissections of $P_ {n+1}$.\vadjust{\goodbreak}

\item[(iii)]It would be interesting to understand what happens when
the sequence
$(\mu_i)_{i \geq2}$ does not satisfy
$(H_\theta)$, for instance, if
$\sum_{i=2}^{\infty} i\mu_i = \infty$. We hope to investigate this in
future work.
\end{longlist}
\end{rems}

\subsection{Galton--Watson trees conditioned on their number of leaves}\label{sec3.1}

Let $\tau\in\mathbb{T}$. Recall our notation $(u(i), 0 \leq i \leq
\zeta(\tau)-1)$ for the vertices of $ \tau$ listed in
lexicographical order
and denote
the number of children of $u(j)$ by $k_j$. Define $\Lambda_{\tau}(l)$
for every $\ell\in\{0,1,\ldots,\zeta(\tau)\}$ by
\[
\Lambda_{\tau}(\ell)=\sum_{0\leq j<\ell}1_{\{k_j=0\}}.
\]
Note that if $Z=W(\tau)$ is the Lukasiewicz path of $\tau$, $\Lambda_{\tau}$
coincides with $\Lambda$ as defined in Proposition \ref{deflien}.
Also note that $\Lambda_{\tau}( \zeta(\tau))=\lambda(\tau)$ is
the total number
of leaves of $\tau$.

\begin{thmm}[(\cite{K})]\label{thmleaves} Let $ ( \mu_j)_ { j \geq2}$
be a sequence of nonnegative real numbers satisfying the assumption
$(H_ \theta)$ for some $ \theta\in(1,2]$. Put $ \mu_1=0$ and $\mu_0=1-\sum_{j=2}^{\infty} \mu_j$,
so that $ \mu= ( \mu_j)_ { j \geq0}$
is a critical probability measure on $ \N$.
For every $n\geq1$ such that $\Pmu[ \lambda(\tau)= n]>0$, let
$\mathfrak
{t}_n$ be a
random tree distributed according to $\Pmu[ \cdot| \lambda(\tau)=
n]$. The following two properties hold:
\begin{longlist}[(ii)]
\item[(i)]We have
\[
\sup_{0 \leq t \leq1} \biggl\llvert \frac{\Lambda_{\mathfrak
{t}_n}(\lfloor\z(\mathfrak{t}_n)
t\rfloor)}{n} - t\biggr\rrvert \mathop{
\longrightarrow}^{( \P)}_ {n
\rightarrow\infty} 0.
\]
\item[(ii)] There exists a sequence $(B_k)_ {k \geq1}$ of positive
constants converging to $ \infty$ such that
%
\begin{equation}
\label{eqcvluka} \biggl(\frac{1}{B_{\z(\mathfrak{t}_n)} } W_{ \lfloor
\z(\mathfrak{t}_n) t \rfloor}(
\mathfrak{t}_n); 0 \leq t \leq1 \biggr) \mathop{\longrightarrow}^{(d)}_{n \rightarrow\infty}
\bigl(X^{\mathrm
{exc}}_t; 0 \leq t \leq 1\bigr).
\end{equation}
\end{longlist}
\end{thmm}

\begin{pf} Note that $\Lambda_ { \mathfrak{t}_n}(\zeta( \mathfrak
{t}_n) )=\lambda(\mathfrak{t}_n)=n$.
In \cite{K}, Corollary 3.3, it is shown that, for every $0 < \eta<1$,
\[
\sup_{\eta\leq t \leq1} \biggl\llvert \frac{\Lambda_ { \mathfrak
{t}_n}(\lfloor\zeta( \mathfrak{t}_n) t\rfloor)}{ \zeta( \mathfrak{t}_n) t} - \mu_0\biggr
\rrvert \mathop {\longrightarrow }^{(\P)}_ { n \rightarrow\infty} 0.
\]
In particular, this implies that $\zeta( \mathfrak{t}_n)/n$ converges in
probability to $1/\mu_0$.
Assertion~(i) follows from the preceding convergences, noting that, for
every $t\in(0,1]$,
\[
\frac{\Lambda_{\mathfrak{t}_n}(\lfloor\z(\mathfrak{t}_n) t\rfloor
)}{n} - t = t \frac{\zeta
(\mathfrak{t}_n)}{n} \biggl( \frac{\Lambda_ { \mathfrak
{t}_n}(\lfloor\zeta( \mathfrak{t}_n) t\rfloor
)}{ \zeta( \mathfrak{t}_n) t} -
\mu_0 \biggr) +t \biggl( \frac{\mu_0 \zeta(\mathfrak{t}_n)}{n} -1 \biggr).
\]
The second assertion is a particular case of
\cite{K}, Theorem 6.1.
\end{pf}

\subsection{Convergence to the stable lamination}\label{sec3.2}

We fix a sequence of nonnegative real numbers $ ( \mu_j)_ { j \geq
2}$ satisfying Assumption $(H_ \theta)$ for some $ \theta\in(1,2]$
and we define $ \mu_0$ and $ \mu_1$ as previously. Throughout this
section, for every $n \geq1$ such that $Z_n$
defined by \eqref{eqzn} is positive (so that $\P^\mu_n$ is well
defined), $\mathfrak l_n$ stands for a random dissection distributed
according to
the Boltzmann probability measure $\P^\mu_n$, and $\mathfrak{t}_n$
stands for its dual tree $\phi(\mathfrak l_n)$, which is
distributed according to $\P_\mu[ \cdot| \lambda(\tau)=n]$ by
Proposition \ref{propGW}. The total progeny of $\mathfrak{t}_n$
is denoted by $\z_n$. The Lukasiewicz path of $\mathfrak{t}_n$ is
denoted by
$W^n$ and $u_0^n,u_1^n, \ldots,u_{\z_n-1}^n$ are the vertices of
$\mathfrak{t}_n$ listed in lexicographical order. Let $(B_n)_ {n \geq
1}$ be a
sequence of positive real numbers such that~\eqref{eqcvluka} holds.
Define the rescaled
Lukasiewicz path $X^n$ of $\mathfrak{t}_n$ by $X^n_t=\frac{1}{B_{\z_n}}
W^n_{\lfloor\z_n t\rfloor}$ for $0 \leq t \leq1$. By Theorem
\ref{thmleaves} and Skorokhod's representation
theorem (see, e.g.,~\cite{Bill}, Theorem~6.7), we may and will assume
that the following convergence holds almost surely in the
space $\R\otimes\D([0,1],\R)$:
%
\begin{equation}
\label{eqCvSk} \biggl( \sup_{0 \leq t \leq1} \biggl\llvert \frac{\Lambda_{\mathfrak{t}_n}(\lfloor\z_n t\rfloor)}{n+1} - t
\biggr\rrvert, X^n \biggr) \,\mathop{\longrightarrow}_{n
\rightarrow\infty}^{a.s.}\,
\bigl(0, X^{\exc} \bigr).
\end{equation}

\subsubsection{Convergence to the Brownian triangulation}

Here, we suppose that \mbox{$\theta=2$}.

\begin{prop}\label{propconvergenceBrownian}When $n$ tends to infinity,
$ D(W^n) \mathop{\rightarrow}\limits^{a.s.} \mathbf{L}(\mathbh{e})$ in the
sense of the Hausdorff distance $d_H$ between compact subsets of
$\Db$.
\end{prop}

\begin{pf}We fix $ \omega$ in the underlying probability space so that
the convergence~\eqref{eqCvSk} holds for this value of $ \omega$ and
we will verify that for this particular value of $ \omega$ we have also
$ D(W^n) \to\mathbf{L}( \mathbh{e})$. Since the space $(\C,d_H)$ is
compact, we may find a random subsequence $(n_k( \omega))$ (depending
on $ \omega$) such that $D(W^ {n_k})$ converges to a compact subset $K$
of $ \Db$, and we need to verify that $K = \mathbf{L}(\mathbh{e})$.
Since $D(W^ {n_k})$ is a dissection for every $k$, one easily checks
that $K$ must be a geodesic lamination of $ \Db$. Since $ \mathbf
{L}(\mathbh{e})$ is a maximal lamination of $ \Db$, the proof will be
complete if we can verify that $ \mathbf{L}(\mathbh{e}) \subset K$.

So we let $0 \leq s < t \leq1$ be such that $s
\stackrel{\mathbh{e}}{\thicksim} t$ and we aim at proving that
$[e^{-2 \textnormal{i} \pi s},e^{-2 \textnormal{i} \pi t}] \subset K$.
Let $ \varepsilon>0$. Simple arguments using the convergence \eqref
{eqCvSk} (and the fact that local minima of the Brownian excursion are
distinct) show that for every $n$ large enough, we can find integers
$i_n, j_n \in\{1, \ldots, \zeta_n-1\}$ such that $ | i_n / \zeta_n-s
| < \varepsilon$, $ | j_n / \zeta_n-t | < \varepsilon$ and
\[
W^n_ {i_n}> W^n_ {i_n-1},\qquad
j_n= \min\bigl\{k> i_n; W^n_k <
W^n_ {i_n}\bigr\}.
\]
By Proposition \ref{propdiscretesons}, $u^n_ {i_n}$ and $u^n_ {j_n}$
are consecutive children of $u^n_ {i_n-1}$. Recalling that
$\Lambda_ { \mathfrak{t}_n}(\zeta( \mathfrak{t}_n) )=n$, we get
from Lemma \ref{lemecriture} that
\[
\biggl[ \exp \biggl(-2 \textnormal{i} \pi\frac{\Lambda_{\mathfrak
{t}_n}(i_n)}{n+1} \biggr), \exp
\biggl(-2 \textnormal{i} \pi \frac{\Lambda_{\mathfrak{t}_n}(j_n)}{n+1} \biggr) \biggr] \subset D
\bigl(W^n\bigr).
\]

To simplify notation, set $s_n = {\Lambda_{\mathfrak
{t}_n}(i_n)}/({n+1}) $ and
$t_n={\Lambda_{\mathfrak{t}_n}(j_n)}/({n+1})$. From the convergence
\eqref
{eqCvSk}, we get $|s_ {n}-s| < \varepsilon$ and $|t_ {n}-t| < \varepsilon$
for every large enough~$n$. In particular, we see that the chord
$[e^{-2 \textnormal{i} \pi s},e^{-2 \textnormal{i} \pi t}]$ lies
within distance $2 \varepsilon$ from $D( W^ {n})$ for every large enough
$n$. It follows that the chord $[e^{-2 \textnormal{i} \pi s},e^{-2
\textnormal{i} \pi t}]$ is within distance $2 \varepsilon$ from $K$. Since
$ \varepsilon>0$ was arbitrary, we get that $[e^{-2 \textnormal{i} \pi s}
,e^{-2 \textnormal{i} \pi t}] \subset K$, which completes the proof.
\end{pf}

\subsubsection{Convergence to the stable lamination when $\theta\neq2$}

We now assume that $ \theta\in(1,2)$. Recall that the convergence
\eqref{eqCvSk} is assumed to hold a.s.

\begin{prop}\label{propconvergence}We have $ D(W^n) \mathop
{\rightarrow}\limits^{a.s.} L(X^\exc)$ as $n \rightarrow\infty$ in the
sense of the Hausdorff distance $d_H$ between compact subsets of
$\Db$.
\end{prop}

We fix $ \omega$ in the underlying probability space so that both the
conclusion of Proposition \ref{propfermL} and the convergence \eqref
{eqCvSk} hold for this value of $ \omega$ and,
furthermore, the path $X^\exc(\omega)$ satisfies properties (H0)--(H4).
We then consider a subsequence $(n_k( \omega))$ such that $D(W^
{n_k})$ converges to a compact subset $K$ of $ \Db$, and we need to
verify that $K = {L}(X^\exc)$. We will first prove that $L(X^\exc
)\subset
K$ before proving the reverse inclusion. In both cases, the precise
description of $L(X^\exc)$ as a union of chords will be crucial. Note
that $K$ must contain the circle $ \mathbb{S}^1$ because the
dissection $D(W^n)$ contains the polygon $P_ {n+1}$. We
stress that the lamination $L(X^\exc)$ is not maximal, in contrast
to the case $\theta=2$. As a consequence, we will have to prove the
nontrivial reverse inclusion.

\begin{lem}\label{lemtechnique}Let $s$ be a jump time of $X^\exc$ and
$t=\inf\{u > s; X^\exc_u = X^\exc_{s-} \}$. For $\varepsilon\in
(0,(t-s)/2)$ small enough, we can choose
an integer $n_0(\varepsilon)$ such that, for every $n \geq n_0(\varepsilon
)$, there
exists $s_n \in(s-\varepsilon,s+\varepsilon) \cap\z_n^ {-1} \N$ such that
the following
inequalities hold:
%
\begin{equation}
\label{eqt3}\inf_{[t-\varepsilon,t+\varepsilon]} X^n < X^n_{s_n -}
< \inf_{[s_n, t-\varepsilon]}X^n.
\end{equation}
\end{lem}

Lemma \ref{lemtechnique} follows from the convergence of $X^n$ to $
X^{\mathrm{exc}}
$ and well-known properties of the Skorokhod topology. We give only the
main ideas of the proof and leave the details to the reader. The time
$s_n$ can be chosen (arbitrarily close to $s$ when~$n$ is large) so
that $X^n_ {s_n-}$ is close to $ X^{\mathrm{exc}}_ {s-}$ and $ \Delta
X^n_ {s_n}$ is
close to $ \Delta X^{\mathrm{exc}}_s$. Then~\eqref{eqt3} is derived
by observing
that, for $ \varepsilon>0$ small enough,
\[
\inf_ {[t,t+ \varepsilon]} X^{\mathrm{exc}}< X^{\mathrm{exc}}_t=
X^{\mathrm{exc}}_ {s-} < \inf_ {[s,t-
\varepsilon]} X^{\mathrm{exc}}.
\]
Notice that the bound $\inf_ {[t,t+ \varepsilon]} X^{\mathrm{exc}}<
X^{\mathrm{exc}}_t$ holds because
otherwise $t$ would be a time of local minimum of $X$
and this would contradict (H4).

\begin{lem}\label{leminclusion1}We have $L(X^\exc) \subset K$.
\end{lem}

\begin{pf} Since $K$ is closed, the property of Proposition \ref
{propfermL} shows that
it is enough to verify that $ [e^{-2 \textnormal{i} \pi\a}, e^{-2
\textnormal{i} \pi\b}  ] \subset K$
for every $(\a,\b)\in\mathcal{E}_1$. So let $(\a,\b)\in\mathcal
{E}_1$. Then $\alpha$ is a jump time of $X^{\mathrm{exc}}$ and
$\b=\inf\{u > \a; X^{\mathrm{exc}}_u = X^{\mathrm{exc}}_{\a-}\}$.
To show that $ [e^{-2
\textnormal{i} \pi\a}, e^{-2 \textnormal{i} \pi
\b}  ] \subset K$, it is sufficient to show that for every $
\varepsilon>0$ and every $n$
sufficiently large we can find $\alpha_ {n},\b_ {n} \in[0,1]$ such that
$|\a_ {n}-\a|\leq2\varepsilon, |\b_ {n}-\b|\leq2\varepsilon$ and $
[e^{-2
\textnormal{i} \pi\a_ {n}}, e^{-2 \textnormal{i} \pi\b_ {n}}
]
\subset D(W^ {n})$. We fix $ \varepsilon>0$. Using Lemma \ref
{lemtechnique} with $(s,t)=( \a, \b)$, we can, for every large enough
$n$, find $\a'_n \in
(\a-\varepsilon,\a+\varepsilon)\cap\zeta_n^{-1}\mathbb{N}$ such that
\[
\inf_{ [\b-\varepsilon,\b+\varepsilon ]}X^n <X^n_{\a'_n-}<
\inf_{ [\a'_n,\b -\varepsilon ]}X^n.
\]
Then put $\b'_n= \inf\{ u \geq\a'_n;
X^n_u<X^n_{\a'_n-} \}$ and note that
$|\alpha-\a'_n| \leq\varepsilon$, $|\b-\b'_n| \leq\varepsilon$. The
time $
\zeta_n\a'_n$
must correspond to a positive jump of $W^n$, and we have also
\[
\zeta_n \b'_n = \inf\bigl\{l\geq
\zeta_n\a'_n ; W^n_l
=W^n_{ \zeta_n
\a'_n} - \bigl(W^n_{ \zeta_n \a'_n}-W^n_{ \zeta_n \a'_n-1}+1
\bigr)\bigr\}.
\]
Using formula (\ref{eqDZ}) and recalling that $\Ln$ coincides with
the process
$\Lambda$ of Proposition \ref{deflien} if $Z=W^n$, we get from Lemma
\ref{lemecriture} that
\[
\biggl[ \exp \biggl(-2 \textnormal{i} \pi\frac{\Ln(\z_n \alpha'_n)}{n+1} \biggr), \exp
\biggl(-2 \textnormal{i} \pi\frac{\Ln(\z_n \beta'_n)}{n+1} \biggr) \biggr] \subset D
\bigl(W^n\bigr).
\]
If we set $\alpha_n=(n+1)^ {-1}\Ln(\z_n \alpha'_n)$
and $\b_n=(n+1)^ {-1} \Ln(\z_n \b'_n)$, the convergence~(\ref
{eqCvSk}) shows that $\a_ {n}$ and $\b_ {n}$ satisfy $|\a_{n}-\a
|\leq
2 \varepsilon$ and
$|\b_{n}- \b| \leq2 \varepsilon$ for all sufficiently large $n$, thus giving
the desired result.
\end{pf}

We now prove the reverse inclusion.

\begin{lem} \label{leminclusion2}We have $K \subset L(X^\exc)$.
\end{lem}

\begin{pf}
Recall that $D(W^ {n_k})$ converges to $K$ in the Hausdorff sense. By
the formula of Proposition \ref{deflien}, we can write
\[
D\bigl(W^ {n_k}\bigr)= \bigcup_ {(u,v) \in\mathcal{E}_ {(n_k)}} \bigl[
e^ {-2
\mathrm{i} \pi u}, e^ {-2 \mathrm{i} \pi v} \bigr],
\]
where $\mathcal{E}_ {(n_k)}$ is a (finite) symmetric subset of
$[0,1]^2$. By extracting a subsequence if necessary, we may assume that
$\mathcal{E}_ {(n_k)} \rightarrow\mathcal{E}_ \infty$
in the Hausdorff sense as $ k \rightarrow\infty$, where $ \mathcal
{E}_ \infty$ is a symmetric closed subset of $[0,1]^2$. It is easy to
verify that
\[
K=\bigcup_ {(u,v) \in\mathcal{E}_ { \infty}} \bigl[ e^ {-2 \mathrm
{i} \pi u},
e^ {-2 \mathrm{i} \pi v} \bigr].
\]
The proof of the inclusion $K \subset L(X^\exc)$ then reduces to
checking that if $ u,v \in\mathcal{E}_ \infty$ with $u<v$, we have $u
\mathop{\simeq}^{X^{\mathrm{exc}}} v$.\vadjust{\goodbreak}

So fix $u,v \in\mathcal{E}_ \infty$ such that $u<v$. Then the pair
$(u,v)$ is the limit of a sequence $(u_k,v_k)$ with
$(u_k,v_k)\in\mathcal{E}_ {(n_k)}$ for every $k$. From Proposition
\ref{deflien}, we can find integers $l_ {n_k}<m_ {n_k}$ in $ \{0,1,
\ldots, \z_ {n_k}\}$ such that
\[
u= \lim_ {k \to\infty} \frac{ \Lambda_{\mathfrak{t}_{n_k}}(l_
{n_k})}{n_k+1},\qquad v=\lim_ {k \to\infty}
\frac{\Lambda_{\mathfrak{t}_{n_k}}(m_ {n_k})}{n_k+1}
\]
and
%
\begin{equation}
\label{eqmnk} m_ {n_k}= \inf \bigl\{ i \geq l_ {n_k};
W^ {n_k}_i=W^
{n_k}_ {l_ {n_k}}-1 \bigr\}.
\end{equation}
By \eqref{eqCvSk}, we have also
%
\begin{equation}
\label{equv} u = \lim_ {k \rightarrow\infty} \frac{l_ {n_k}}{ \z_
{n_k}},\qquad v =
\lim_ {k \rightarrow\infty} \frac{m_ {n_k}}{ \z_ {n_k}}.
\end{equation}
From \eqref{eqmnk}, we have $W^ {n_k}_i \geq W^ {n_k}_ {m_ {n_k}}$ for
every $ i \in[l_ {n_k}, m_ {n_k}]$. Thus, using the convergence of
$X^n$ to $ X^{\mathrm{exc}}$ and \eqref{equv},
%
\begin{equation}
\label{eqxs} X^{\mathrm{exc}}_s \geq X^{\mathrm{exc}}_ {v-}\qquad
\mathrm{for\ every}\ s \in(u,v).
\end{equation}
From property (H3) this implies that $ X^{\mathrm{exc}}_v = X^{\mathrm
{exc}}_ {v-}$, and then
$(B_{\zeta_{n_k}})^{-1}W^{n_k}_{m_{n_k}}= X^{n_k}_{m_{n_k}/ \zeta_{n_k}}$
must converge to $X^{\mathrm{exc}}_v$. Note that $X^{\mathrm
{exc}}_{u-}$ and $X^{\mathrm{exc}}_u$ are the only
possible accumulation points for the sequence
$(B_{\zeta_{n_k}})^{-1}W^{n_k}_{l_{n_k}}= X^{n_k}_{l_{n_k}/ \zeta_{n_k}}$.
Now consider two cases:

\begin{longlist}[$-$]
\item[$-$]
If $ X^{\mathrm{exc}}_u= X^{\mathrm{exc}}_ {u-}$, then $(B_{\zeta
_{n_k}})^{-1}W^{n_k}_{l_{n_k}}=
X^{n_k}_{l_{n_k}/ \zeta_{n_k}}$ converges
to $X^{\mathrm{exc}}_u$ and, using~\eqref{eqmnk}, we get that $ X^{\mathrm{exc}}_u= X^{\mathrm
{exc}}_v$. It follows that $ X^{\mathrm{exc}}_s >
X^{\mathrm{exc}}_v$ for every $s \in(u,v)$, because otherwise this would contradict
(H1) or (H4). Clearly, we obtain $u \mathop{\simeq}^{X^{\mathrm
{exc}}} v$.

\item[$-$]
If $ X^{\mathrm{exc}}_u > X^{\mathrm{exc}}_ {u-}$, then we must have
$X^ {n_k}_ {l_ {n_k}/ \z_
{n_k}} \rightarrow X^{\mathrm{exc}}_ {u-}$
[otherwise \eqref{eqmnk} would give $X^{\mathrm{exc}}_u=X^{\mathrm
{exc}}_v$, and \eqref{eqxs}
would contradict (H2)]. Then \eqref{eqmnk} gives $ X^{\mathrm
{exc}}_v= X^{\mathrm{exc}}_ {u-}$.
The inequality \eqref{eqxs} can then be reinforced in $ X^{\mathrm
{exc}}_s > X^{\mathrm{exc}}_v=
X^{\mathrm{exc}}_ {u-}$ for every $s \in(u,v)$, since otherwise $X^{\mathrm{exc}}$
would have a
local minimum equal
to $X^{\mathrm{exc}}_v=X^{\mathrm{exc}}_{u-}$ in $(u,v)$, which would
contradict (H4).
Hence, we also get $u \mathop{\simeq}^{X^{\mathrm{exc}}} v$ in that case.
\end{longlist}
This completes the proof.
\end{pf}

Together with Lemmas \ref{leminclusion1}, \ref{leminclusion2}
completes the proof of Theorem \ref{thmcv1} in the case $ \theta\neq2$.

\subsection{\texorpdfstring{Description of the faces of $L(X^{\mathrm{exc}})$ for $\theta\neq2$}
{Description of the faces of L(X exc) for theta not equal 2}}\label{sec3.3}

We still consider the case $1< \theta< 2$. By definition, the faces of
$L(X^{\mathrm{exc}})$ are the connected components of $\Db\setminus
L(X^{\mathrm{exc}})$.
In this section, we study the faces of $L(X^{\mathrm{exc}})$ and we
show in
particular that, almost surely, every face of $L(X^{\mathrm{exc}})$ is
bounded by
infinitely many chords (in contrast to the case $ \theta=2$ where all
faces are triangles).

\begin{lem}\label{lemcomponent}Almost surely, for every face $U$ of
$L(X^{\mathrm{exc}})$, if $\Gamma=\mathbb{S}_1 \cap\overline{U} $
denotes the part
of the
boundary of $U$ lying on the circle, then:
\begin{longlist}[(iii)]
\item[(i)] $U$ is a convex open set;
\item[(ii)] $\Gamma$ is not a singleton;
\item[(iii)] $1\notin\Gamma$.
\end{longlist}
\end{lem}
\begin{pf}
Assertions (i) and (ii) hold for any geodesic lamination of $ \Db$, and
we leave the proof to the reader. To get (iii), fix $ \varepsilon>0$ and
note that by Proposition
\ref{propapproxX} we can find $s \in(0, \varepsilon]$ and $t \in[1-
\varepsilon, 1)$ such that the chord $[e^{-2 \textnormal{i} \pi s},e^{-2
\textnormal{i} \pi t}]$ is contained in $L(X^{\mathrm{exc}})$. It
follows that $1$
cannot belong to the boundary of a connected component of $\Db
\setminus
L(X^{\mathrm{exc}})$.
\end{pf}

For distinct
$s,t \in(0,1)$, we denote by $\mathbb{H}^s_t$ the open half-plane
bounded by the line containing $e^{-2 \textnormal{i} \pi s}$ and $e^{-2
\textnormal{i} \pi t}$
and such that $1 \notin\mathbb{H}^s_t$. We write $\widetilde
{\mathbb
{H}}^s_t$ for the other open half-plane bounded
by the same line.

\begin{prop}\label{propcomponents}Let $s$ be a jump time of
$X^{\mathrm{exc}}$ and
$t=\inf\{u > s; X^{\mathrm{exc}}_u =
X^{\mathrm{exc}}_{s-}\}$. There exists a unique face $U$ of
$L(X^{\mathrm{exc}})$ contained in
$\mathbb{H}^s_t$ and whose closure
$\overline{U}$ contains the chord $[e^{-2 \textnormal{i} \pi s},e^{-2
\textnormal{i} \pi t}]$. The face $U$ is called the face associated to
$s$. The mapping $s \mapsto U$ is a one-to-one correspondence between
jump times of $ X^{\mathrm{exc}}$ and faces of $L(X^{\mathrm{exc}})$.
\end{prop}

\begin{pf}
We start by giving a
description of the face associated to $s$. Let $(\alpha_i,\beta_i)_{i
\geq1}$ be defined by
\begin{eqnarray*}
\bigl\{(\alpha_i,\beta_i); i \geq1\bigr\}&=& \Bigl\{(\a,
\b) ; s \leq\a< \b\leq t,
 X_\a=X_\b=
\inf_ {[s, \a]} X \mbox{ and }
 X^{\mathrm{exc}}_r > X^{\mathrm{exc}}_\a\\
 &&\hspace*{202pt}\mbox{for } r\in(\a,\b)\Bigr\},
\end{eqnarray*}
where the pairs $(\alpha_i,\beta_i)$ are listed in such a way that
$\beta_i-\alpha_i> \beta_j-\alpha_j$ for
$i<j$. The intervals $(\alpha_i,\beta_i)$ are exactly the excursion
intervals of $(X_{r}-I^s_{r})_{s\leq r\leq t}$ away from $0$.
Note that $ \a_i
\mathop{\simeq}^{X^{\mathrm{exc}}} \b_i$ by Proposition \ref
{propidentifications},
and that the intervals $(\alpha_i,\beta_i)$, $i\geq1$
are disjoint. Furthermore, the fact that (H3) holds for $X^{\mathrm{exc}}$
shows that the times $\alpha_i$, $i\geq1$ are not jump times of
$X^{\mathrm{exc}}$.

For every $n\geq1$, let $V_n$ be the convex open polygon whose vertices
are
\[
\bigl\{e^{-2\textnormal{i} \pi s},e^{-2\textnormal{i} \pi t}\bigr\} \cup \bigcup
_{i=1}^n \bigl\{e^{-2\textnormal{i} \pi\a_i},e^{-2\textnormal{i} \pi\b_i}
\bigr\}.
\]
Observe that $V_n\subset V_{n+1}$. We finally set
\[
V=\bigcup_{n \geq1} V_n,
\]
which is a convex open set. It is clear that $V$ is
contained in the open half-plane $\mathbb{H}^s_t$ and that
$\overline{V}$ contains $[e^{-2\textnormal{i} \pi s},e^{-2\textnormal
{i} \pi t}]$. To prove that
$V$ is a connected component of $\Db\setminus L(X^{\mathrm{exc}})$,
we proceed
in two steps. We first prove that $V \subset\Db\setminus
L(X^{\mathrm{exc}})$
and then that $V$ is a maximal connected open subset of $\Db
\setminus L(X^{\mathrm{exc}})$.

Let us prove that $V \subset\Db\setminus
L(X^{\mathrm{exc}})$. Argue by contradiction and suppose that there
exist $P \in
L(X^{\mathrm{exc}})$ and $N \geq1$ such that $P \in V_N$. By the
definition of
$L(X^{\mathrm{exc}})$, there exist $0 \leq u \leq v < 1$ such that $u
\mathop{\simeq}^{X^{\mathrm{exc}}} v$ and $P \in[e^{-2\textnormal
{i} \pi u},
e^{-2\textnormal{i} \pi v}]$.
Since $V$ is contained in the open half-plane $\mathbb{H}^s_t$, we
must have $s \leq u < v \leq t$. Let us first show that $s<u$. If
$s=u$, the definition of $\mathop{\simeq}^{X^{\mathrm{exc}}}$
implies that
$v=\inf\{r > s; X^{\mathrm{exc}}_r = X^{\mathrm{exc}}_{s-}\}=t$.
Consequently, $P \in
[e^{-2\textnormal{i} \pi s},e^{-2\textnormal{i} \pi t}]$, contradicting
the fact that $P \in
V_N$. We thus have $s<u$. Since $P \in V_N$ and since for every $j \in
\{1, \ldots, N\}$ the chord $  [ e^{-2\textnormal{i} \pi\a
_j},e^{-2\textnormal{i} \pi\b_j}  ]$ does not cross the chord $
[ e^{-2\textnormal{i} \pi u},e^{-2\textnormal{i} \pi v}  ]$,
a simple argument shows that there exists $1 \leq i
\leq N$ such that
$u \leq\a_i < \b_i \leq v$, the case $(u,v)=(\a_i,\b_i)$ being
excluded. 
We examine two cases:
\begin{longlist}[$-$]
\item[$-$] If $u< \a_i$, then $X^{\mathrm{exc}}_{u-} > X^{\mathrm
{exc}}_{\a_i}$ because $ \inf_{[s,\a
_i]} X^{\mathrm{exc}}= X^{\mathrm{exc}}_{\a_i}$, $\a_i$ is
a local minimum time for $X^{\mathrm{exc}}$ and local minima are
almost surely
distinct. Since $\a_i \in[u,v]$ and $u \mathop{\simeq}^{X^{\mathrm
{exc}}} v$,
this contradicts Remark \ref{remeq}.
\item[$-$] If $u= \a_i$, we know that $u$ is not a jump time of $
X^{\mathrm{exc}}$ and
the property $u \mathop{\simeq}^{X^{\mathrm{exc}}} v$ implies $v =
\inf\{r>u; X^{\mathrm{exc}}_r
\leq X^{\mathrm{exc}}_ { \a_i}\}= \b_i$, which is excluded.
\end{longlist}
In each case, a contradiction occurs. This completes the first step.

Let us then prove that $V$ is a maximal connected open
subset of $\Db\setminus L(X^{\mathrm{exc}})$. To this end, we
observe that we have
\[
V = \mathbb{H}^s_t \cap \Biggl( \bigcap
_{i=1}^\infty\widetilde {\mathbb {H}}^ { \a_i}_ { \b_i}
\Biggr) \cap\D.
\]
The fact that $V$ is contained in the set in the right-hand side is
immediate from our construction, and the
reverse inclusion is also easy. Set $R= ( \mathbb{H}^s_t)^ {c} \cap
\Db
$ and $R_i= (\widetilde{ \mathbb{H}}^ { \a_i}_ { \b_i})^ {c} \cap
\Db$
for $i \geq1$.
It follows that
%
\begin{equation}
\label{complement-face} \Db\setminus V= \mathbb{S}_1 \cup R \cup
\Biggl( \bigcup_ {i
=1}^\infty R_i
\Biggr).
\end{equation}
This implies that the boundary of $V$ is contained in $L( X^{\mathrm
{exc}})$, and it
follows that $V$ is a maximal connected open
subset of $\Db\setminus L(X^{\mathrm{exc}})$. From the preceding
formula for $\Db
\setminus V$, it is also clear that
the boundary of $V$ contains the chord $[e^{-2 \textnormal{i} \pi
s},e^{-2 \textnormal{i} \pi t}]$, as well as all chords
$[e^{-2 \textnormal{i} \pi\a_i},e^{-2 \textnormal{i} \pi b_i}]$, and
we have obtained the existence of the
face associated to $s$. The uniqueness of this face is obvious for
geometric reasons.

We still have to prove the last assertion of the proposition.
Let $U$ be a face of $L(X^{\mathrm{exc}})$. We need to verify that $U$
is the face
associated to a certain
jump time of $X^{\mathrm{exc}}$. To this end,\vadjust{\goodbreak} let
$\Gamma=\mathbb{S}_1 \cap\overline{U}$ be the part of the boundary
of $U$ lying on the circle and set:
\[
s= \inf\bigl\{ u \geq0 ; e^{-2 \textnormal{i} \pi u} \in\Gamma\bigr\},\qquad t=\sup\bigl\{ 0 \leq
u \leq1 ; e^{-2 \textnormal{i} \pi u} \in\Gamma\bigr\}.
\]
By Lemma \ref{lemcomponent}(iii), we have $0<s < t < 1$. By the
compactness of $L( X^{\mathrm{exc}})$ and a convexity argument, it is easy
to verify that $ [e^{-2 \textnormal{i} \pi s}, e^{-2 \textnormal
{i} \pi t} ] \subset
L(X^{\mathrm{exc}})$. We then claim that $s$ is a jump time of
$X^{\mathrm{exc}}$. If not, by
Proposition \ref{propidentifications}, this means that $X^{\mathrm
{exc}}_s=X^{\mathrm{exc}}_t$
and $X^{\mathrm{exc}}_u>X^{\mathrm{exc}}_s$ for $u \in(s,t)$. But
then Proposition
\ref{propapproxX} could be used to produce a chord of $L( X^{\mathrm{exc}})$
partitioning $U$ into two
disjoint open sets, which is impossible. So $s$ is a jump time of $
X^{\mathrm{exc}}$
and we then know
that $t=\inf\{u > s; X^{\mathrm{exc}}_u
=X^{\mathrm{exc}}_{s-}\}$. Let $V$ be the face associated to
$s$. To prove that $U=V$, it is sufficient to show that $U
\cap V \neq\varnothing$. This follows from simple
geometric considerations.
This completes the proof.
\end{pf}

\section{The stable lamination coded by a continuous function}\label{sec4}

The definitions of the limiting random laminations
$\mathbf{L}(\mathbh{e})$ and $L(X^{\mathrm{exc}})$ that appear in
our main result
Theorem \ref{thmcv1}
for $\theta=2$ and $\theta\neq2$ were somewhat different. The goal
of this
section is to unify these two cases by explaining how $L(X^{\mathrm
{exc}})$ (for
$\theta\neq2$) can also be constructed from a random continuous function.
This will allow us to make the connection between our stable laminations
and the so-called stable trees, which were studied in particular in
\cite{DuquesneLG,DuquesneLG-fractal},
in the same way as the Brownian triangulation is connected to the
Brownian CRT
\cite{Aldous2}, and this will also be useful when we calculate the
Hausdorff dimension of $L(X^{\mathrm{exc}})$. The relevant random
function, called
the height process in continuous time, was introduced in
\cite{LeJan} and studied in great detail in
\cite{DuquesneLG}.

In this section, $X$ is the strictly stable spectrally positive
L\'{e}vy of index $\theta$, as defined in Section \ref{sec2.1} and $1<\theta<2$.

\subsection{The height process}\label{sec4.1}

The continuous-time height process associated with~$X$ can be defined
by the following approximation formula. For every $ t \geq0$,
\[
H_t= \lim_ { \varepsilon\rightarrow0}\frac{1}{\varepsilon} \int
_0^t ds \mathbh{1}_{\{X_s \leq I^s_t+\varepsilon\}},
\]
where the convergence holds in probability. The process $(H_t)_ {t \geq
0}$ has a continuous modification, which we consider from now on.

A very useful ingredient in the study of the height process is the
so-called exploration process $ ( \rho_t)_ { t \geq0}$, which is a
strong Markov process taking values in the space $M_f ( \R_+)$ of all
finite measures on $ \R_+$. For every $ t \geq0$, $ \rho_t$ is
defined by
%
\begin{equation}
\label{eqrho}\langle\rho_t, f\rangle= \int_{[0,t]}
d_s I^s_t f({H_s})
\end{equation}
for every measurable function $f\dvtx \R_+ \rightarrow\R_+$. Here the
notation $d_s I^s_t$ refers to the integration\vadjust{\goodbreak} with respect to the
nondecreasing function $s \rightarrow I^s_t$ (recall the definition of
$I^s_t$ in Section \ref{sec2.1}). Note, in particular, that $ \langle\rho_t,1 \rangle=X_t-I_t$. The process $ ( \rho_t)_ {t \geq0}$ enjoys the
following two important properties \cite{DuquesneLG}, Lemma~1.2.2:
\begin{longlist}[(ii)]
\item[(i)] Almost surely for every $t \geq0$, $\rho_t( \{0\})=0$ and $
\operatorname{supp}( \rho_t)=[0,H_t]$ [here and later $\operatorname{supp}(
\mu)$ denotes the topological support of $ \mu\in M_f ( \R_+)$, with
the convention that $ \operatorname{supp}(0)= \{0\}$].
\item[(ii)] Almost surely $ \{t \geq0; H_t=0\}= \{ t \geq0; \rho_t=0\}
= \{t \geq0; X_t=I_t\}$.
\end{longlist}
In addition to (i), one can prove that, for every fixed $t \geq0$, $
\rho_t( \{H_t\})=0$ almost surely. This follows from formula (17) in
\cite{DuquesneLG}. Moreover, almost surely for every jump time $s$ of
$X$, $ \rho_s( \{H_s\})= \Delta X_s$ (see formula (19) in \cite{DuquesneLG}).

We will need another important property of the exploration process. To
state this property, we need to introduce some notation. If $ \mu\in
M_f ( \R_+)$ and $ \alpha\geq0$, the ``killed'' measure $ k_ \a\mu$
is the unique element of $M_f( \R_+)$ such that, for every $t \geq0$,
\[
k_ \a\mu\bigl( [0,t]\bigr)= \mu\bigl([0,t]\bigr) \wedge\bigl( \mu(
\R_+) - \alpha\bigr)^+.
\]
Suppose that $ \mu\in M_f ( \R_+)$ has compact support and set $h(
\mu
)= \sup( \operatorname{supp}( \mu))$. Then if $ \nu\in M_f( \R_+)$,
the concatenation $ [ \mu, \nu] \in M_f ( \R_+)$ is defined by
\[
\bigl\langle[ \mu, \nu],f\bigr\rangle= \langle\mu,f \rangle+ \int\nu (dt) f
\bigl(h( \mu)+t\bigr).
\]

Let $T$ be a stopping time of the filtration of $X$ and let $X_t^
{(T)}=X_ {T+t}-X_T$ for every $t \geq0$.
Recall that $(X^{(T)}_t)_{t\geq0}$ has the same distribution as
$(X_t)_{t\geq0}$ by the
strong
Markov property of $X$. Set $I^ {(T)}_t= \inf_ {s \leq t}X^ {(T)}_s$
for every $t \geq0$, and let $(H_t^ {(T)})_ {t \geq0}$ and $ (\rho^
{(T)}_t)_ {t \geq0}$ be, respectively, the height process and the
exploration process associated with $X^ {(T)}$. According to formula
(20) in \cite{DuquesneLG}, we have almost surely for every $t \geq0$,
%
\begin{equation}
\label{eqr1} \rho_ {T+t}= \bigl[ k_ {- I^ {(T)}_t}
\rho_T, \rho_t^
{(T)} \bigr].
\end{equation}
It follows that almost surely for every $t \geq0$,
%
\begin{equation}
\label{eqr2}H_ {T+t}- \inf_ {s \in[T,T+t]} H_s=
H_t^ {(T)}
\end{equation}
(see \cite{DuquesneLG}, Lemma 1.4.5, for the case where $T$ is
deterministic, but the derivation is the same in the general case).

The following result is a continuous analog of Proposition
\ref{propdiscretesons}.

\begin{prop}\label{propptsmultiples}The following holds almost surely.
Let $s \geq0$ be a jump
time of $X$ and $t=\inf\{u > s; X_u = X_{s-}\}$. Then:
\begin{longlist}[(iii)]
\item[(i)] for every $u \in[s,t]$, $H_u \geq H_s$ and $H_u=H_s$ if and
only if
$X_u=\inf_{[s,u]}X$;
\item[(ii)]for every $\a\in[0,s)$, $\inf_{[\a,s]}H< H_s$;
\item[(iii)] for every $u \in(t,\infty)$, $\inf_{[s,u]}H< H_s$.\vadjust{\goodbreak}
\end{longlist}
\end{prop}

\begin{pf}
Since the set of all jump times can be written as a countable
collection of stopping times, it is sufficient to consider the case
when $s=S$ is a stopping time, that is, also a jump time of $X$, and $t=
T= \inf\{r \geq S; X_r= X_ {S-}\}$. By preceding observations, we
know that $ \rho_S ( \{H_S\})= \Delta X_S$.

Let us prove (i). From \eqref{eqr1} applied to the stopping time $S$,
we have $ \rho_ {S+r} \geq k_ { \Delta X_S} \rho_S$ for every $r \in
[0, T-S]$ and, thus,
\[
H_ {S+r}= \sup( \operatorname{supp } \rho_ {S+r}) \geq\sup(
\operatorname {supp } k_ { \Delta X_S} \rho_S)=H_S.
\]
Furthermore, for the same values of $r$, \eqref{eqr1} shows that $H_
{S+r}=H_S$ can only hold if $ \rho^ {(S)}_r=0$, which is equivalent [by
\eqref{eqrho}] to $X^ {(S)}_r=I^ {(S)}_r$. This completes the proof of (i).

To get (ii), we observe that we can always pick a rational $ \beta\in
( \alpha, S)$ such that $X_ \beta< X_S$. By \eqref{eqr2} applied to
$T= \beta$,
\[
H_S - \inf_ {r \in[ \a,S]} H_r \geq
H_S-\inf_ {r \in[ \b,S]} H_r= H^
{( \beta)}_ {S- \beta}.
\]
Since $X_S > X_ \b$, we have $ \langle\rho^ {( \b)}_ {S- \b},1
\rangle
\geq X^ {( \beta)}_ {S- \beta}>0$ and, thus, $H^ { ( \beta)}_ {S-
\beta
}>0$, completing the proof of (ii).

Finally, for every $ \varepsilon>0$ set $T_ \varepsilon= \inf\{r \geq S;
X_r \leq
X_ {S-}- \varepsilon\}$. By \eqref{eqr1} we have $ \rho_ {T_ \varepsilon
}=k_ { \Delta
X_s + \varepsilon} \rho_S$ and $H_ {T_ \varepsilon}= \sup( \operatorname
{supp } k_ {
\Delta X_s+ \varepsilon} \rho_S)<H_S$ because $ \rho_S( \{H_S\})=
\Delta
X_S$. This completes the proof.
\end{pf}

The following result will also be useful.

\begin{prop} \label{propClasseH2} The following holds almost surely
for every choice of
$0 \leq s <t$ such that $H_s=H_t$ and $H_u>H_s$ for all
$u \in(s,t)$. For every $ \varepsilon\in(0,(t-s)/2)$, there
exist $ s' \in(s,s+\varepsilon)$ and $ t'\in(t-\varepsilon,t)$ such that
$s'< t'$ and:
\begin{longlist}[(ii)]
\item[(i)] $H$ does not attain a local minimum at $s'$ or at
$t'$;
\item[(ii)] $H_{s'}=H_{t'}=\inf_{[s',t']} H$ and there exists $v \in
(s',t')$ such that $H_v=H_{s'}$.
\end{longlist}
\end{prop}

\begin{pf}
We can assume that $\varepsilon< (t-s)/4$. Set
$m=\inf_{[s+\varepsilon,t-\varepsilon]}H$. By the continuity of $H$,
there exists $\varepsilon' \in(0, \varepsilon)$ such that
$\sup_{[s,s+\varepsilon']}H < m$. Let $u \in(s,s+\varepsilon') \cap\Q$.
We have
\[
\inf_{ [u, s +
\varepsilon']} H< H_u
\]
because it easily follows from formula \eqref{eqr1} that $ \inf_
{[q,q+ \delta]} H < H_q$ for every rational $q>0$ and every $ \delta
>0$, almost surely (the point is that the measure $ \rho_q$ gives no
mass to $ \{H_q\}$, so that the supremum of the support of $ k_a \rho_q$ will be strictly smaller than $H_q$, for every $a>0$).

Then let $v \in(u, s + \varepsilon']$ be such that
$H_v= \inf_{[u, s+ \varepsilon']} H$. Finally, set $s'=\inf\{ r \in[s,s
+ \varepsilon']; H_r=H_v\}$ and $t' = \sup\{ r \in[s + \varepsilon', t]
; H_r=H_v\}$ so that~$H$ does not attain a local minimum at $s'$
or at $t'$. By construction and using the continuity of~$H$,
we have
\[
s < s'\leq u<v \leq s+\varepsilon< t-\varepsilon< t'
< t.
\]
Since $H_{s'}=H_v=H_{t'}$, the proposition is proved.
\end{pf}

\subsection{The normalized excursion of the height process}\label{sec4.2}

Recall the notation of Section \ref{sec2.1}, where we
have constructed the normalized excursion $X^{\mathrm{exc}}$
from the excursion of $X$ straddling $1$.

The normalized excursion of the height process, which is denoted by
$H^\exc$, is defined as follows. Set $ \beta_ \varepsilon= \theta/(
\Gamma(2-
\theta) \varepsilon^ { \theta-1})$. Using
Proposition \ref{proplaw}, one shows that there exists a continuous process
$(H^{\mathrm{exc}}_t)_{0 \leq t \leq1}$, such that, for every $t$
belonging to a
subset of $[0,1]$
of full Lebesgue measure,
\[
H_t^\exc=\lim_{\varepsilon\rightarrow0}\frac{1}{\b_{\varepsilon}} \card
\Bigl\{ u \in[0,t]; X^{\mathrm{exc}}_{u-} < \inf_{ [u,t
]}X^{\mathrm{exc}},
\Delta X^{\mathrm{exc}}_u > \varepsilon\Bigr\} \qquad\mbox{a.s.}
\]
See \cite{Duquesne}, Section 3, for details of the argument. This
process $H^{\mathrm{exc}}$ is called the
normalized excursion of the height process. The pair $(X^{\mathrm
{exc}},H^\exc)$ can
be constructed explicitly from the process $X$
via the formula
%
\begin{equation}
\label{defHexc} \quad \bigl(X^{\mathrm{exc}}_t,H^\exc_t
\bigr)_{0\leq t\leq1} = \bigl(\zeta_1^{-{1}/{\theta}}
(X_{\underline{g}_1+\zeta_1 t} - X_{\underline{g}_1}), \zeta_1^{({1}/{\theta})-1}
H_{\underline{g}_1+\zeta_1 t} \bigr)_{0\leq t\leq1},
\end{equation}
where we recall the notation $\underline{g}_1=\sup\{s \leq1; X_s=I_s
\}$ and $\zeta_1=\underline{g}_1- \inf\{s > 1; X_s=I_s\}$.

\begin{rem}\label{remform} From formula \eqref{defHexc}, we see that
the results of Propositions \ref{propptsmultiples} and \ref{propClasseH2}
remain valid if we replace $X$ with $X^{\mathrm{exc}}$ and $H$ with
$H^\exc$. More
precisely, we will use these results in the following
form. Almost surely:
\begin{enumerate}[(1)]
\item[(1)] Let $0 \leq s \leq1$ be a jump
time of $X^{\mathrm{exc}}$ and $t=\inf\{u > s; X^{\mathrm{exc}}_u =
X^{\mathrm{exc}}_{s-}\}$. Then for $u
\in[s,t]$, $ H^{\mathrm{exc}}_u \geq H^{\mathrm{exc}}_s$, and
$H^{\mathrm{exc}}_u=H^{\mathrm{exc}}_s$ if and only if $X^{\mathrm{exc}}_u=\inf_{[s,u]}X^{\mathrm{exc}}$.
Moreover, if $0 \leq\a<s$, then $\inf_{[\a,s]}H^{\mathrm{exc}}<
H^{\mathrm{exc}}_s$, and if
$t<u\leq1$, then $\inf_{[s,u]}H^{\mathrm{exc}}< H^{\mathrm{exc}}_s$;
\item[(2)] For every choice of $0 \leq s<t \leq1$, the conditions
$H^{\mathrm{exc}}_s=H^{\mathrm{exc}}_t$ and $H^{\mathrm{exc}}_u>H^{\mathrm{exc}}_s$
for all $u \in(s,t)$
imply that for every $\varepsilon>0$ sufficiently small, there exist $
s' \in(s,s+\varepsilon)$ and $ t'\in(t-\varepsilon,t)$ such that:
\begin{enumerate}[(ii)]
\item[(i)] $H^{\mathrm{exc}}$ does not attain a local minimum at $s'$
or at
$t'$,
\item[(ii)] $\inf_{[s',t']} H^{\mathrm{exc}}= H^{\mathrm
{exc}}_{s'}=H^{\mathrm{exc}}_{t'}$ and there exists $u
\in(s',t')$ such that $H^{\mathrm{exc}}_u=H^{\mathrm
{exc}}_{s'}=H^{\mathrm{exc}}_{t'}$.
\end{enumerate}
\end{enumerate}
\end{rem}

The main result of \cite{Duquesne} states that if $\mathfrak{t}_n$ is
a $GW_\mu$ tree conditioned on having
total progeny $n$, the discrete height process $  (H_{k}(\mathfrak
{t}_n)  )_ {0 \leq k \leq n}$,
appropriately rescaled, converges in distribution to $H^{\mathrm
{exc}}$. However,
we will not use this fact.



\subsection{Laminations coded by continuous functions}\label{sec4.3}

Let $g\dvtx [0,1] \rightarrow\R_+$ be a continuous function such that
$g(0)=g(1)=0$. We define a
pseudo-distance on $[0,1]$ by
\[
d_g(s,t)=g(s)-g(t)- 2 \min_{r \in[s
\wedge t, s\vee t]}g(r)
\]
for $s,t \in[0,1]$. The associated
equivalence relation on $[0,1]$ is defined by setting $s
\stackrel{g}{\thicksim} t$ if and only if $d_g(s,t)=0$ or,
equivalently, $g(s)=g(t)=\break \min_{r \in[s \wedge t, s\vee
t]}g(r)$ (in the special case $g=\mathbh{e}$, this equivalence
relation was already used in
Section \ref{sec2}).

The quotient set $T_g := [0,1] / \stackrel{g}{\thicksim}$ equipped with
the distance $d_g$
is an $\R$-tree, called the tree coded by the function $g$. We refer to
\cite{DuquesneLG-fractal,Evans}
for more information about
$\R$-trees, which are natural generalizations of discrete trees, and
their coding by functions.

For $s \in[0,1]$, we let $\cl_g(s)$ be the equivalence class of $s$
with respect to the equivalence relation
$ \stackrel{g}{\thicksim} $. Then, for $s,t \in[0,1]$, we set $
s \stackrel{g}{\thickapprox} t$ if at least one of the following two
conditions holds:
\begin{longlist}[$-$]
\item[$-$] $ s \mathop{\thicksim}\limits^g t$ and $g(r)>g(s)$ for every $r
\in
(s\wedge t, s\vee t)$;
\item[$-$] $ s \mathop{\thicksim}\limits^g t$ and $s\wedge t= \min\cl_g(s)$,
$s\vee
t=\max\cl_g(s).$
\end{longlist}

By \cite{CLG}, Proposition 2.5, the set
\[
\mathbf{L}(g):= \bigcup_{s \stackrel{g}{\thickapprox} t}\bigl[e^{-2
\textnormal{i} \pi s}
,e^{-2 \textnormal{i} \pi t}\bigr]
\]
is a geodesic lamination of $\Db$. Note that if $g=\mathbh{e}$, this
coincides with the definition in Section \ref{sec2}, thanks to the fact
that local minima of $\mathbh{e}$ are distinct.

In what follows we take $g=H^{\mathrm{exc}}$ and write $\thickapprox
^{H^{\mathrm{exc}}}$ rather
than $ \stackrel{H^{\mathrm{exc}}}{\thickapprox}$ for notational reasons.

\begin{prop}\label{propinfini} Almost surely, for every real number $u
\in[0,1]$ such that $\card( \cl_{H^{\mathrm
{exc}}}(u)) \geq
3$, there exists a jump time $\a$ of $X^{\mathrm{exc}}$ such that $\a
\in
\cl_{H^{\mathrm{exc}}}(u)$. Conversely, let $\a$ be a
jump time of
$X^{\mathrm{exc}}$
and $\b= \inf\{
r > \a; X^{\mathrm{exc}}_r=X^{\mathrm{exc}}_{\a-}\}$. Then
$\card(\cl_{H^{\mathrm{exc}}}(\a)) =
\infty$, furthermore,
$\min\cl_{H^{\mathrm{exc}}}(\a)=\a$ and $\max
\cl_{H^{\mathrm{exc}}}(\a
)$ $=\b$,
so that, in particular, $\a
\thickapprox^{H^\exc} \b$.
\end{prop}

\begin{pf} The first assertion is a consequence of Theorem 4.7 in \cite
{DuquesneLG-fractal} and the discussion following
this statement. The fact that $\card(\cl_{H^{\mathrm{exc}}}(\a)) =
\infty$ if $\a$ is a jump time of $X^{\mathrm{exc}}$ follows from
\cite{DuquesneLG-fractal}, Theorem
4.6. Finally, let $\a$ be a jump time of
$X^{\mathrm{exc}}$ and
let $\b=\inf\{ r \geq\a;
X^{\mathrm{exc}}_r=X^{\mathrm{exc}}_{\a-}\}$. By the first part of
Remark \ref{remform}, we
know that
$H^{\mathrm{exc}}_\a=\inf_{ [\a,\b ]}H^{\mathrm
{exc}}=H^{\mathrm{exc}}_\b$ and that for any
$\varepsilon>0$,\looseness=-1
\[
\inf_{ [\a-\varepsilon,\a ]} H^{\mathrm{exc}}< H^{\mathrm
{exc}}_\a,\qquad
\inf_{ [\b,\b+\varepsilon ]} H^{\mathrm{exc}}< H^{\mathrm{exc}}_\b.
\]\looseness=0
The desired result follows.\vadjust{\goodbreak}
\end{pf}

\begin{thmm}\label{thmXH}Almost surely, the relations
$\mathop{\simeq}^{X^\exc} $ and $\mathop{\thickapprox}^{H^{\mathrm{exc}}}$ coincide. In particular,
\[
L\bigl(X^{\mathrm{exc}}\bigr)=\mathbf{L}\bigl(H^{\mathrm{exc}}\bigr)
\qquad\mbox{a.s.}
\]
\end{thmm}

\begin{pf} We first observe that both relations $\mathop{\simeq}^{X^\exc} $ and $\mathop{\thickapprox}^{H^{\mathrm{exc}}}$
are closed, in the sense that their graphs are
closed subsets of
$[0,1]^2$. In the case of
$\mathop{\simeq}^{X^\exc} $, this was already observed in the proof of
Proposition \ref{proplamination}.
In the case of $\mathop{\thickapprox}^{H^{\mathrm{exc}}}$, this is elementary (see \cite{CLG}, Section 2.3).

Let $s,t\in[0,1]$ such that $s<t$ and $s \mathop{\simeq}^{X^\exc} t$.
From Proposition
\ref{propapproxX}, we can write the pair $(s,t)$ as the limit of a
sequence $(s_n,t_n)$ in $\mathcal{E}_1$
(of course, if $s$ is a jump time of $X^{\mathrm{exc}}$, we take
$s_n=s$ and $t_n=t$
for every $n$). However,
Proposition~\ref{propinfini} then implies that
$s_n
\thickapprox^{H^{\mathrm{exc}}} t_n$, for every $n$, and it follows
that $s \thickapprox^{H^{\mathrm{exc}}} t$.

Let us prove the converse. Let $(s,t)$ be such that $0\leq s<t\leq1$
and $s \mathop{\thickapprox}^{H^{\mathrm{exc}}} t$. If $\card(\cl_{H^{\mathrm{exc}}}(s)) \geq
3$, we must have $s=\min\card(\cl_{H^{\mathrm{exc}}}(s))$ and
$t= \max\card(\cl_{H^{\mathrm{exc}}
}(s))$, so that
Proposition \ref{propinfini} implies that the pair $(s,t)$ belongs to
$\mathcal{E}_1$, and, in
particular, $s \mathop{\simeq}^{X^\exc} t$. If $\card(\cl_{H^{\mathrm{exc}}
}(s)) =
2$, then the second part of Remark \ref{remform} shows that
$(s,t)$ is the limit of a sequence of pairs $s_n,t_n$ such that
$s_n \mathop{\thickapprox}^{H^{\mathrm{exc}}} t_n$ and $\card(\cl_{H^{\mathrm{exc}}}(s_n))
\geq
3$. We have then $s_n \mathop{\simeq}^{X^\exc} t_n$ for every $n$ and
$s \mathop{\simeq}^{X^\exc} t$ since the relation $\mathop{\simeq}^{X^\exc} $ is closed.
\end{pf}

\begin{rem} In the discrete setting, the definition of the dissection
$D(W(\tau))$ via formula
\eqref{eqDZ} uses the times $s^i_1,\ldots,s^i_{k_i}$, which
can be defined either from the Lukasiewicz path of $\tau$ as in
Proposition \ref{propdiscretesons}(i) or from the discrete height
process of $\tau$ as in part (ii)
of the same proposition. In the continuous setting, we recover these
two different points of view in the definition of the $\theta$-stable
lamination as $L(X^{\mathrm{exc}})$ or $\mathbf{L}(H^{\mathrm{exc}})$.
\end{rem}

\section{The Hausdorff dimension of the stable lamination}\label{sec5}

In this section we determine the Hausdorff dimension of $L(X^{\mathrm
{exc}})$ and of
some other random
sets related to $L(X^{\mathrm{exc}})$.
We refer the reader to \cite{Mattila} for background concerning
Hausdorff and Minkowski dimensions.

\begin{thmm}\label{thmdimH}Fix $\theta\in(1,2]$. Let $L(X^{\mathrm
{exc}})$ be the
random lamination coded by the normalized
excursion $X^{\mathrm{exc}}$ of the $\theta$-stable L\'{e}vy process
and let $A$
stand for the set of all endpoints of chords in $L(X^\exc)$. Then
\[
\dim(A) = 1 -\frac{1}{\theta},\qquad \dim\bigl( L\bigl(X^{\mathrm
{exc}}\bigr)
\bigr)=2-\frac{1}{\theta},
\]
where $\dim(K)$ stands for the Hausdorff dimension of a subset $K$ of
$\mathbb{C}$.
Furthermore, if $1<\theta<2$, then a.s. for every face $V$ of
$L(X^{\mathrm{exc}})$,
\[
\dim\bigl(\overline V\cap\mathbb{S}^1\bigr)= \frac{1}{\theta}.
\]
\end{thmm}

\begin{rem}In the case $\theta=2$, the results of the theorem are
already known; See \cite{Aldous}
for a sketch of the argument and
\cite{LGP} for a detailed proof. We thus restrict our attention to
$\theta\in
(1,2)$. We follow the
idea of the proof of \cite{LGP} but a different argument is needed
because of the existence of jump times.
\end{rem}

It will be convenient to identify the interval $[0,1)$ with $\mathbb
{S}_1$ via the mapping $x
\mapsto e^{-2\textnormal{i} \pi x}$. The set $A$ of the theorem is the
set of all $x \in
\mathbb{S}^1$ such that there exists $y \in\mathbb{S}^1$ with $y
\neq x$
and $x \mathop{\simeq}^{X^\exc} y$. We also let
$\mathcal{I}$ be the set of all (ordered) pairs $(I,J)$, where $I$
and $J$ are two disjoint closed subarcs of $\mathbb{S}_1$ with
nonempty interior and rational endpoints. If $(I,J) \in
\mathcal{I}$, we denote by $A^{(I,J)}$ the set of all $x \in I$ such
that $x \mathop{\simeq}^{X^\exc} y$ for some $y \in J$. In
particular,
\[
A=\bigcup_{(I,J) \in\mathcal{I}} A^{(I,J)}.
\]

In the following, $\underline{\dim}_M(B)$ and $\overline{\dim}_M(B)$ will denote, respectively, the lower and the upper Minkowski
dimensions of a set $B$ (see \cite{Mattila} for definitions). In order
to compute Hausdorff and Minkowski dimensions, the following
proposition will be useful.

\begin{prop}\label{proprange}Almost surely, for every $t > 0$, the set
$\{ 0 \leq s \leq t ;
S_s=X_s\}$ has Hausdorff dimension and upper Minkowski dimension
equal to $1-1/\theta$, and the set $\{0\leq s\leq t; I_s=X_s\}$
has Hausdorff dimension and upper Minkowski dimension
equal to $1/\theta$.
\end{prop}

\begin{pf} Recall that if $(\tau_t, t \geq0)$ is a stable subordinator
of parameter
$\rho\in(0,1)$, then, almost surely, for all $t>0$, the Hausdorff
dimension and the upper Minkowski dimension of $\{ \tau_s;
0 \leq s \leq t\}$, or of the closure of this set, is equal to $\rho$
(see, e.g.,
\cite{BertoinSub}, Theorem 5.1, Corollary 5.3). Let $L=(L_t, t \geq0)$
stand for a local time of $S-X$
at $0$, and let $L^ {-1}$ be the right-continuous inverse of $L$.
Since $X$ has only positive jumps, the set $\{ 0 \leq s < t ; S_s=X_s\}
$ is
closed. By \cite{Bertoin}, Lemma VIII.1, $L^{-1}$ is a subordinator
of index $1-1/\theta$ and by \cite{Bertoin}, Proposition~IV.7, $\{ 0
\leq
s < t ; S_s=X_s\}$ coincides with the closure of $\{ L^{-1}_s; 0
\leq s < L_t\}$. As $L_t >0$ almost surely, the first assertion of the
proposition follows.
The proof of the second assertion is similar, noting that $-I$ is a
local time at $0$
for $X-I$ and that the right-continuous inverse of $-I$ is a
stable subordinator of index $1/\theta$, again by \cite{Bertoin}, Lemma
VIII.1.
\end{pf}

\begin{lem}\label{lemdim}For $a \in(0,1]$, set $\widehat{F}_a:=\{ u
\in(0,a); X^{\mathrm{exc}}_{u-} \leq\inf_{ [u,a ]}
X^{\mathrm{exc}}\}$. Almost surely, for every jump
time $a$ of $ X^{\mathrm{exc}}$ in $(0,1)$ we have
%
\begin{equation}
\label{eqdim}\dim(\widehat{F}_a)= \overline{\dim }_M(
\widehat{F}_a)=1-\frac{1}{\theta}.
\end{equation}
\end{lem}

Informally, if one identifies the interval $[0,1]$ with the circle $
\mathbb{S}_1$ by using the map $ x \rightarrow e^ {-2 \mathrm{i} \pi
x}$, the set
$\widehat{F}_a$ corresponds to endpoints in $(0,a)$ of chords that
connect a point
of $(0,a)$ to a point of $(a,1)$.

\begin{pf*}{Proof of Lemma \ref{lemdim}} We first consider an analog
of $\widehat{F}_a$ where $X^{\mathrm{exc}}$ is replaced by the
L\'evy process $X$. Precisely, for every $a>0$, we set
\[
\widetilde{F}_a:=\Bigl\{u\in(0,a); X_{u-}\leq
\inf_{[u,a]}X\Bigr\}.
\]
Note that, under the condition $X_a>I_a$, $\widetilde{F}_a$ is
contained in the (closure of the)
excursion interval of $X-I$ that straddles $a$. Thanks to this
observation and to the connection between
$X^{\mathrm{exc}}$ and $X$ given by Proposition \ref{proplaw}, the
result of the lemma
will follow if we can verify that
%
\begin{equation}
\label{eqdimbis}\dim(\widetilde{F}_a)= \overline {\dim
}_M(\widetilde{F}_a)=1-\frac{1}{\theta}
\end{equation}
for every jump time $a$ of $X$ [note that if $ X^{\mathrm{exc}}$ is
given by the
formula of Proposition~\ref{proplaw}, the jump times of $ X^{\mathrm
{exc}}$ exactly
correspond to jump times of $X$ over $( \underline{g}_1, \underline
{d}_1)$]. Let $K>0$ and consider only jump times that are bounded above
by~$K$. The desired result for such jump times follows by considering
the process $X$ time-reversed at time $K$ and using the strong Markov
property together with Proposition~\ref{proprange}.
\end{pf*}

\begin{pf*}{Proof of Theorem \ref{thmdimH}} We first prove the last
assertion of the theorem. By Proposition
\ref{propcomponents}, a face $V$ of $L(X^{\mathrm{exc}})$ is
associated to a jump
time $s$ of $X^{\mathrm{exc}}$, and we set
$t=\inf\{r>s\dvtx X^{\mathrm{exc}}_{r}=X^{\mathrm{exc}}_{s-}\}$. Let the
intervals $(\alpha_i,\beta_i)$, $i\geq1$ be defined as in the proof
of Proposition~\ref{propcomponents}. Then, it easily follows from
\eqref{complement-face} that
\[
\overline V \cap\mathbb{S}^1 = [s,t]\Bigm\backslash \bigcup
_{i=1}^\infty (\alpha_i,
\beta_i) = \Bigl\{r\in[s,t]; X^{\mathrm{exc}}_r=
\inf_{[s,r]} X^{\mathrm{exc}}\Bigr\},
\]
where we recall that $\mathbb{S}^1$ is identified with $[0,1)$.
The calculation of $\dim(\overline V \cap\mathbb{S}^1)$ now follows
from the second assertion of Proposition \ref{proprange},
using also Proposition~\ref{proplaw}.

Let us turn to the first part of the theorem. We follow the ideas of
the proof
of the analogous result in \cite{LGP}. We will prove that
%
\begin{equation}
\label{eqdim2}\dim(A)=1-1/\theta, \qquad\overline{\dim}_M
\bigl(A^{(I,J)}\bigr) \leq1-1/\theta
\end{equation}
for every $(I,J) \in\mathcal{I}$, a.s. If (\ref{eqdim2})
holds, then
\begin{eqnarray*}
\underline{\dim}_M\bigl(A^{(I,J)} \cup A^{(J,I)}
\bigr) &\leq& \overline{\dim}_M\bigl(A^{(I,J)} \cup
A^{(J,I)}\bigr)
\\
&=&\max\bigl(\overline{\dim}_M\bigl(A^{(I,J)}\bigr),
\overline{\dim}_M\bigl(A^{(J,I)}\bigr)\bigr)
\\
&\leq& \dim(A),
\end{eqnarray*}
and then the same argument as in Proposition 2.3 of \cite{LGP} entails
that
\[
\dim\bigl(L\bigl( X^{\mathrm{exc}}\bigr)\bigr)=1+ \dim(A)=2-1/ \theta.
\]
It remains to establish
(\ref{eqdim2}). In order to verify that
\[
\overline{\dim}_M \bigl(A^{(I,J)}\bigr) \leq1 - 1/\theta
\]
for every $(I,J) \in\mathcal{I}$,
we need only consider the case $I=[u,v], J=[u',v']$ with $0\leq u' <
v' \leq1$, $0\leq u < v \leq1$ (if one of the subarcs $I$ or $J$
contains $0$ as an interior point, partition it into two subarcs
whose interior does not contain $0$). Since the relations
$\mathop{\simeq}^{X^\exc}$ and $\thickapprox^{H^\exc}$
coincide, the time-reversal invariance property of
$H^{\mathrm{exc}}$ (see \cite{DuquesneLG}, Corollary 3.1.6) allows us
to restrict
to the case $0 \leq u < v < u' < v' \leq1$. Choose a jump time $a$
of $X^{\mathrm{exc}}$ such that $v < a < u'$ and observe that
$\widehat{F}_a
\subset A$ and $A^{(I,J)} \subset\widehat{F}_a$, with the notation
of Lemma \ref{lemdim}. Hence, by the latter lemma,
$\overline{\dim}_M (A^{(I,J)}) \leq\overline{\dim}_M
(\widehat{F}_a) = 1-1/\theta$. Lemma \ref{lemdim} and the property
$\widehat{F}_a \subset A$ also give $ 1-1/\theta\leq\dim A$. We
have then
\[
1-\frac{1}{\theta} \leq\dim A \leq \overline{\dim}_M (A) \leq
\max_{(I,J) \in\mathcal{I}} \overline{\dim}_M \bigl(A^{(I,J)}\bigr)
\leq1 - \frac{1}{\theta}.
\]
In particular, $\dim A=1-1/\theta$ and \eqref{eqdim2} holds. This
completes the proof.
\end{pf*}

\section*{Acknowledgments}
I am deeply indebted to Jean-Fran\c{c}ois Le
Gall for suggesting me to study this model, for insightful
discussions and for carefully reading the manuscript and making many
useful suggestions.

%


\printaddresses

\end{document}